\documentclass{autart}
\usepackage{natbib}

\usepackage{amsmath,amssymb}
\usepackage{pstricks,pst-plot}
\usepackage{graphicx,subfigure,xspace,bm}
  
\newtheorem{theorem}{Theorem}[section]

\newtheorem{lemma}[theorem]{Lemma}

\newtheorem{remark}[theorem]{Remark}
\newtheorem{example}[theorem]{Example}

\newtheorem{proposition}[theorem]{Proposition}

\newcommand{\wdin}{\upsubscr{d}{in}{w}}

\newcommand{\wdout}{\upsubscr{d}{out}{w}}

\newcommand{\Nin}{\upscr{\mathcal{N}}{in}}
\newcommand{\Nout}{\upscr{\mathcal{N}}{out}}

\newcommand{\WW}{\mathcal{W}}
\renewcommand{\SS}{\mathcal{S}_{\text{agree}}}
\newcommand{\smallnonzeroeig}{\Lambda_{*}}

\newcommand{\realpart}{\mathrm{Re}}
\newcommand{\impart}{\mathrm{Im}}

\newcommand{\real}{{\mathbb{R}}}
\newcommand{\realpositive}{{\mathbb{R}}_{>0}}
\newcommand{\realnonnegative}{{\mathbb{R}}_{\ge 0}}
\newcommand{\integers}{\mathbb{Z}}
\newcommand{\integerspositive}{\mathbb{Z}_{\geq 1}}

\newcommand{\fdist}{\mathfrak{f}}

\newcommand{\until}[1]{\{1,\dots,#1\}}
\newcommand{\map}[3]{#1:#2 \rightarrow #3}
\newcommand{\setmap}[3]{#1:#2 \rightrightarrows #3}
\newcommand{\equilibria}[1]{\operatorname{Eq}(#1)}

\newcommand{\network}{\Sigma}

\newcommand{\networkeng}{\subscr{\Sigma}{eng}}

\newcommand{\Bgraph}{\mathcal{G}}

\newcommand{\directedNash}{\subscr{\Psi}{Nash-dir}}

\newcommand{\Lie}{\mathcal{L}}
\newcommand{\SetLie}{\widetilde{\mathcal{L}}}
\newcommand{\gradient}{\nabla}
\newcommand{\Bconvexhull}[1]{\mathrm{co}\Big\{ #1\Big\}}
\newcommand{\ones}{\mathbf{1}}
\newcommand{\zeros}{\mathbf{0}}
\newcommand{\setdef}[2]{\{#1 \; | \; #2\}}

\newcommand{\Lyap}{V}

\newcommand{\game}[4]{\subscr{\mathbf{G}}{adv-net}=(#1,#2,#3,#4)}
\newcommand{\gametilde}[4]{\subscr{\mathbf{\tilde{G}}}{adv-net}=(#1,#2,#3,#4)}
\newcommand{\gamef}{\gametilde{\network_1}{\network_2}{\networkeng}{\ULift}}
\newcommand{\gameforiginal}{\game{\network_1}{\network_2}{\networkeng}{U}}
\newcommand{\gamew}{\subscr{\mathbf{\tilde{G}}}{adv-net}}
\newcommand{\gameworiginal}{\subscr{\mathbf{G}}{adv-net}}

\newcommand{\undirectedNash}{\subscr{\Psi}{Nash-undir}}

\newcommand{\Blap}{\mathbf{L}}

\newcommand{\fLift}{\tilde{f}}
\newcommand{\ULift}{\tilde{U}}
\newcommand{\gLift}{\tilde{g}}
\newcommand{\xnet}{\bm{x}}
\newcommand{\ynet}{\bm{y}}
\newcommand{\znet}{\bm{z}}

\newcommand{\identity}[1]{\mathsf{I}_{#1}}
\newcommand{\setX}{\mathsf{X}}

\newcommand{\partialvector}{\partial}
\newcommand{\Dout}{\subscr{\mathsf{D}}{out}}

\newcommand{\Adj}{\mathsf{A}}
\newcommand{\Lap}{\mathsf{L}}
\newcommand{\vertices}{\mathcal{V}}
\newcommand{\edges}{\mathcal{E}}

\newcommand\subscr[2]{#1_{\textup{#2}}}
\newcommand\upscr[2]{#1^{\textup{#2}}}
\newcommand\upsubscr[3]{#1_{\textup{#2}}^{\textup{#3}}}

\newcommand{\oprocendsymbol}{\hbox{$\bullet$}}
\newcommand{\oprocend}{\relax\ifmmode\else\unskip\hfill\fi\oprocendsymbol}

\newcommand{\longthmtitle}[1]{\mbox{}\textup{\textbf{(#1):}}}

\renewcommand{\theenumi}{(\roman{enumi})}

\newcommand{\myclearpage}{\clearpage}
\renewcommand{\myclearpage}{}

\begin{document}

\runauthor{B. Gharesifard and J. Cort\'es}

\begin{frontmatter}

    
  \title{Distributed convergence to Nash equilibria in
    \\
    two-network zero-sum games\thanksref{footnoteinfo}}

  \thanks[footnoteinfo]{Incomplete versions of this paper were
    presented as~\citep{BG-JC:11-acc} in the American Control
    Conference~2012 and as~\citep{BG-JC:12-cdc2} in the IEEE Control
    and Decision Conference~2012. This work was performed while
    B. Gharesifard was a postdoctoral researcher at the University of
    California, San Diego.}

  \author[address1]{B. Gharesifard} \ead{bgharesi@illinois.edu}\qquad
  \author[address2]{J. Cort\'es} \ead{cortes@ucsd.edu}
  
  \address[address1]{Coordinated Science Laboratory, University
    of Illinois at Urbana-Champaign, IL, 61801, USA}
  
  \address[address2]{Department of Mechanical and Aerospace Engineering,
    University of California, San Diego, CA, 92093, USA}

  \begin{abstract}
    This paper considers a class of strategic scenarios in which two
    networks of agents have opposing objectives with regards to the
    optimization of a common objective function.  In the resulting
    zero-sum game, individual agents collaborate with neighbors in
    their respective network and have only partial knowledge of the
    state of the agents in the other network.  For the case when the
    interaction topology of each network is undirected, we synthesize
    a distributed saddle-point strategy and establish its convergence
    to the Nash equilibrium for the class of strictly concave-convex
    and locally Lipschitz objective functions.  We also show that this
    dynamics does not converge in general if the topologies are
    directed.  This justifies the introduction, in the directed case,
    of a generalization of this distributed dynamics which we show
    converges to the Nash equilibrium for the class of strictly
    concave-convex differentiable functions with globally Lipschitz
    gradients. The technical approach combines tools from algebraic
    graph theory, nonsmooth analysis, set-valued dynamical systems,
    and game theory.
  \end{abstract}
  
  
\end{frontmatter}

\myclearpage
\section{Introduction}\label{section:intro}

Recent years have seen an increasing interest on networked strategic
scenarios where agents may cooperate or compete with each other
towards the achievement of some objective, interact across different
layers, have access to limited information, and are subject to
evolving interaction topologies. This paper is a contribution to this
body of work.  Specifically, we consider a class of strategic
scenarios in which two networks of agents are involved in a zero-sum
game. We assume that the objective function can be decomposed as a sum
of concave-convex functions and that the networks have opposing
objectives regarding its optimization. Agents collaborate with the
neighbors in their own network and have partial information about the
state of the agents in the other network.  Such scenarios are
challenging because information is spread across the agents and
possibly multiple layers, and networks, by themselves, are not the
decision makers.  Our aim is to design a distributed coordination
algorithm that can be used by the agents to converge to the Nash
equilibrium. Note that, for a 2-player zero-sum game of the type
considered here, a pure Nash equilibrium corresponds to a saddle point
of the objective function.

\textit{Literature review.}  Multiple scenarios involving networked
systems and intelligent adversaries in sensor networks, filtering,
finance, and wireless communications~\citep{JSK-SB:08,PW-MDL:09} can
be cast into the strategic framework described above.  In such
scenarios, the network objective arises as a result of the aggregation
of agent-to-agent adversarial interactions regarding a common goal,
and information is naturally distributed among the agents.  The
present work has connections with the literature on distributed
optimization and zero-sum games.  The distributed optimization of a
sum of convex functions has been intensively studied in recent years,
see e.g.~\citep{AN-AO:09,PW-MDL:09,BJ-MR-MJ:09,MZ-SM:12}. These works
build on consensus-based
dynamics~\citep{ROS-JAF-RMM:07,WR-RWB:08,FB-JC-SM:08cor,MM-ME:10} to
find the solutions of the optimization problem in a variety of
scenarios and are designed in discrete time.  Exceptions
include~\citep{JW-NE:10, JW-NE:11} on continuous-time distributed
optimization on undirected networks and~\citep{BG-JC:12-tac} on
directed networks.

Regarding zero-sum games, the works~\citep{KA-LH-HU:58,DM:77,AS-AO:09}
study the convergence of discrete-time subgradient dynamics to a
saddle point.  Continuous-time best-response dynamics for zero-sum
games converges to the set of Nash equilibria for both
convex-concave~\citep{JH-SS:06} and
quasiconvex-quasiconcave~\citep{EN-RG-RRJ:10} functions.  Under strict
convexity-concavity assumptions, continuous-time subgradient flow
dynamics converges to a saddle
point~\citep{KA-LH:51,KA-LH-HU:58}. Asymptotic convergence is also
guaranteed when the Hessian of the objective function is positive
definite in one argument and the function is linear in the
other~\citep{KA-LH-HU:58,DF-FP:10}.  The distributed computation of
Nash equilibria in noncooperative games, where all players are
adversarial, has been investigated under different assumptions.  The
algorithm in~\citep{SL-TB:87} relies on all-to-all communication and
does not require players to know each other's payoff functions (which
must be strongly convex). In~\citep{PF-MK-TB:12,MSS-KHJ-DMS:12},
players are unaware of their own payoff functions but have access to
the payoff value of an action once it has been executed. These works
design distributed strategies based on extremum seeking techniques to
seek the set of Nash equilibria.


\textit{Statement of contributions.}  We introduce the problem of
distributed convergence to Nash equilibria for two networks engaged in
a strategic scenario.  The networks aim to either maximize or minimize
a common objective function which can be written as a sum of
concave-convex functions.  Individual agents collaborate with
neighbors in their respective network and have partial knowledge of
the state of the agents in the other one.  Our first contribution is
the introduction of an aggregate objective function for each network
which depends on the interaction topology through its Laplacian and
the characterization of a family of points with a saddle property for
the pair of functions. We show the correspondence between these points
and the Nash equilibria of the overall game. When the graphs
describing the interaction topologies within each network are
undirected, the gradients of these aggregate objective functions are
distributed. Building on this observation, our second contribution is
the synthesis of a consensus-based saddle-point strategy for
adversarial networks with undirected topologies. We show that the
proposed dynamics is guaranteed to asymptotically converge to the Nash
equilibrium for the class of strictly concave-convex and locally
Lipschitz objective functions. Our third contribution focuses on the
directed case. We show that the transcription of the saddle-point
dynamics to directed topologies fails to converge in general. This
leads us to propose a generalization of the dynamics, for strongly
connected weight-balanced topologies, that incorporates a design
parameter.  We show that, by appropriately choosing this parameter,
the new dynamics asymptotically converges to the Nash equilibrium for
the class of strictly concave-convex differentiable objective
functions with globally Lipschitz gradients.  The technical approach
employs notions and results from algebraic graph theory, nonsmooth and
convex analysis, set-valued dynamical systems, and game theory.  As an
intermediate result in our proof strategy for the directed case, we
provide a generalization of the known characterization of cocoercivity
of concave functions to concave-convex functions.

The results of this paper can be understood as a generalization to
competing networks of the results we obtained in~\citep{BG-JC:12-tac}
for distributed optimization.  This generalization is nontrivial
because the payoff functions associated to individual agents now also
depend on information obtained from the opposing network.  This
feature gives rise to a hierarchy of saddle-point dynamics whose
analysis is technically challenging and requires, among other things,
a reformulation of the problem as a constrained zero-sum game, a
careful understanding of the coupling between the dynamics of both
networks, and the generalization of the notion of cocoercivity to
concave-convex functions.

\textit{Organization.}  Section~\ref{section:prelim} contains
preliminaries on nonsmooth analysis, set-valued dynamical systems,
graph theory, and game theory.  In
Section~\ref{sec:problem-statement}, we introduce the zero-sum game
for two adversarial networks involved in a strategic scenario and
introduce two novel aggregate objective functions.
Section~\ref{sec:nash-undirected} presents our algorithm design and
analysis for distributed convergence to Nash equilibrium when the
network topologies are undirected.  Section~\ref{sec:nash-directed}
presents our treatment for the directed case.
Section~\ref{sec:conclusions} gathers our conclusions and ideas for
future work.  Appendix~\ref{sec:appendix} contains the generalization
to concave-convex functions of the characterization of cocoercivity of
concave functions.

\myclearpage
\section{Preliminaries}\label{section:prelim}

We start with some notational conventions.  Let $\real$,
$\realnonnegative$, $\integers$, $\integerspositive$ denote the set of
real, nonnegative real, integer, and positive integer numbers,
respectively.  We denote by $ ||\cdot || $ the Euclidean norm on $
\real^d$, $ d \in \integerspositive$ and also use the short-hand
notation $ \ones_d=(1,\ldots,1)^T $ and $\zeros_d=(0,\ldots,0)^T \in
\mathbb{R}^d$.  We let $ \identity{d} $ denote the identity matrix in
$ \mathbb{R}^{d\times d}$.  For matrices $ A\in \real^{d_1\times d_2}
$ and $ B \in \real^{e_1\times e_2} $, $ d_1,d_2,e_1,e_2 \in
\integerspositive $, we let $ A\otimes B $ denote their Kronecker
product.  The function $\map{f}{\setX_1\times \setX_2}{\mathbb{R}} $,
with $ \setX_1 \subset \real^{d_1} $, $ \setX_2\subset \real^{d_2} $
closed and convex, is \emph{concave-convex} if it is concave in its
first argument and convex in the second one~\citep{RTR:97}.  A point $
(x_1^*,x_2^*) \in \setX_1\times \setX_2 $ is a \emph{saddle point} of
$ f $ if $ f(x_1,x_2^*)\leq f(x^*_1,x^*_2)\leq f(x_1^*,x_2)$ for all $
x_1 \in \setX_1 $ and $ x_2 \in \setX_2$. Finally, a set-valued map
$\setmap{f}{\real^d}{\real^d}$ takes elements of $\real^d$ to subsets
of $\real^d$.

\subsection{Nonsmooth analysis}

We recall some notions from nonsmooth analysis~\citep{FHC:83}.  A
function $ \map{f}{\real^d}{\real} $ is \emph{locally Lipschitz} at $
x \in \real^d $ if there exists a neighborhood $ \mathcal{U} $ of $ x
$ and $ C_x \in \realnonnegative$ such that $ |f(y)-f(z)|\leq C_x
||y-z|| $, for $ y,z \in \mathcal{U} $. $f$ is locally Lipschitz on
$\real^d$ if it is locally Lipschitz at $x$ for all $x \in \real^d$
and \emph{globally Lipschitz} on $ \real^d $ if for all $ y,z \in
\real^d $ there exists $ C \in \realnonnegative $ such that $
|f(y)-f(z)|\leq C ||y-z|| $.  Locally Lipschitz functions are
differentiable almost everywhere. The \emph{generalized gradient} of $
f $ is
\[
\partial f(x) = \Bconvexhull{\lim_{k \rightarrow \infty} \nabla f(x_k)
  \ | \ x_k \rightarrow x, x_k \notin \Omega_f \cup S},
\]
where $ \Omega_f $ is the set of points where $ f $ fails to be
differentiable and $ S $ is any set of measure zero.

\begin{lemma}\longthmtitle{Continuity of the generalized gradient
    map}\label{le:gradient_properties}
  Let $ \map{f}{\real^d}{\real} $ be a locally Lipschitz function at $
  x \in \real^d $.  Then the set-valued map $ \setmap{\partial
    f}{\real^d}{\real^d} $ is upper semicontinuous and locally bounded
  at $ x \in \real^d $ and moreover, $ \partial f(x) $ is nonempty,
  compact, and convex.
\end{lemma}

For $\map{f}{\real^d \times \real^d}{\real}$ and $z \in \real^d$, we
let $ \partial_x f(x,z) $ denote the generalized gradient of $x
\mapsto f(x,z)$. Similarly, for $x \in \real^d$, we let $ \partial_z
f(x,z) $ denote the generalized gradient of $z \mapsto f(x,z)$. A
point $ x \in \real^d $ with $ \zeros \in \partial f(x) $ is a
\emph{critical point} of $ f $.  A function $ \map{f}{\real^d}{\real}
$ is \emph{regular} at $ x \in \real $ if for all $ v \in \real^d $
the right directional derivative of $ f $, in the direction of $ v $,
exists at $ x $ and coincides with the generalized directional
derivative of $ f $ at $ x$ in the direction of~$v$. We refer the
reader to~\citep{FHC:83} for definitions of these notions. A convex
and locally Lipschitz function at $ x $ is
regular~\citep[Proposition~2.3.6]{FHC:83}.  The notion of regularity
plays an important role when considering sums of Lipschitz functions.

\begin{lemma}\longthmtitle{Finite sum of
    locally Lipschitz functions}\label{le:finite_sum_gen}
  Let $ \{f^i\}_{i=1}^n$ be locally Lipschitz at $ x \in \real^d
  $. Then $\partial(\sum_{i=1}^nf^i)(x)\subseteq \sum_{i=1}^n\partial
  f^i(x)$, and equality holds if $ f^i $ is regular for $ i \in
  \{1,\ldots, n\} $.
\end{lemma}

A locally Lipschitz and convex function $f$ satisfies, for all $x,x'
\in \real^d$ and $ \xi \in \partial f(x) $, the \emph{first-order
  condition} of convexity,
\begin{equation}\label{eq:convex_prop}
  f(x')-f(x) \geq \xi\cdot (x'-x).
\end{equation}


\subsection{Set-valued dynamical systems}

Here, we recall some background on set-valued dynamical systems
following~\cite{
  JC:08-csm-yo}.  A continuous-time set-valued dynamical system on
$\setX \subset \mathbb{R}^d$ is a differential inclusion
\begin{equation}\label{eq:diff-inc}
  \dot{x}(t) \in \Psi(x(t))
\end{equation}
where $ t \in \realnonnegative $ and $\setmap{\Psi}{\setX \subset
  \real^d}{\real^d}$ is a set-valued map.  A solution to this
dynamical system is an absolutely continuous curve $ x:
[0,T]\rightarrow \setX $ which satisfies~\eqref{eq:diff-inc} almost
everywhere. The set of equilibria of~\eqref{eq:diff-inc} is denoted by
$\equilibria{\Psi} = \setdef{x \in \setX}{ 0 \in \Psi(x)}$.

\begin{lemma}\longthmtitle{Existence of solutions}\label{le:solution}
  For $\setmap{\Psi}{\real^d}{\real^d}$ upper semicontinuous with
  nonempty, compact, and convex values, there exists a solution
  to~\eqref{eq:diff-inc} from any initial condition.
\end{lemma}

The LaSalle Invariance Principle for set-valued continuous-time
systems is helpful to establish the asymptotic stability properties of
systems of the form~\eqref{eq:diff-inc}.  A set $ W \subset \setX $ is
\emph{weakly positively invariant} with respect to $ \Psi $ if for any
$ x \in W $, there exists $ \tilde{x} \in \setX $ such that $
\tilde{x} \in \Psi(x) $.  The set $ W $ is \emph{strongly positively
  invariant} with respect to $ \Psi $ if $ \Psi(x) \subset W $, for
all $ x \in W $.
Finally, the \emph{set-valued Lie derivative} of a differentiable
function $ \map{V}{\real^d}{\real} $ with respect to $ \Psi $ at $ x
\in \real^d $ is defined by $ \SetLie_{\Psi}{V(x)}=\{v \cdot \nabla
V(x) \ | \ v \in \Psi(x) \}$.


\begin{theorem}\longthmtitle{Set-valued LaSalle Invariance
    Principle}\label{th:laSalle}
  Let $ W \subset \setX $ be a strongly positively invariant
  under~\eqref{eq:diff-inc} and $ \Lyap: \setX \rightarrow \mathbb{R}
  $ a continuously differentiable function.  Suppose the evolutions
  of~\eqref{eq:diff-inc} are bounded and $ \max
  \SetLie_{\Psi}{\Lyap(x)} \leq 0 $ or $
  \SetLie_{\Psi}{\Lyap(x)}=\emptyset $, for all $ x \in W $. If
  $ S_{\Psi,\Lyap} = \{x\in \setX \ | \ 0 \in
  \SetLie_{\Psi}{\Lyap(x)}\}, $
  then any solution $ x(t) $, $ t\in \mathbb{R}_{\geq 0} $, starting
  in $W$ converges to the largest weakly positively invariant set $ M
  $ contained in $ \bar{S}_{\Psi,\Lyap}\cap W$. When $ M $ is a finite
  collection of points, then the limit of each solution 
  equals one of them.
\end{theorem}

\subsection{Graph theory}

We present some basic notions from algebraic graph theory following
the exposition in~\citep{FB-JC-SM:08cor}.  A \emph{directed graph}, or
simply \emph{digraph}, is a pair $\Bgraph=(\vertices,\edges)$, where
$\vertices$ is a finite set called the vertex set and $ \edges
\subseteq \vertices\times \vertices $ is the edge set.
A digraph is \emph{undirected} if $(v,u) \in \edges$ anytime $(u,v) \in
\edges$. We refer to an undirected digraph as a \emph{graph}.
A path is an ordered sequence of vertices such that any ordered pair
of vertices appearing consecutively is an edge of the digraph.  A
digraph is \emph{strongly connected} if there is a path between any
pair of distinct vertices. For a graph, we refer to this notion simply
as \emph{connected}.  A \emph{weighted digraph} is a triplet $
\Bgraph=(\vertices,\edges,\Adj) $, where $ (\vertices,\edges) $ is a
digraph and $ \Adj \in \mathbb{R}^{n\times n}_{\geq0} $ is the
\emph{adjacency matrix} of $\Bgraph$, with the property that $
a_{ij}>0 $ if $ (v_i,v_j)\in \edges $ and $ a_{ij}=0 $, otherwise.
The weighted out-degree and in-degree of $v_i$, $i \in \{1,\dots,n\}$,
are respectively, $ \wdout(v_i) =\sum_{j=1}^{n}a_{ij} $ and
$\wdin(v_i)=\sum_{j=1}^n a_{ji} $.  The \emph{weighted out-degree
  matrix} $ \Dout$ is the diagonal matrix defined by $
(\Dout)_{ii}=\wdout(i) $, for all $ i \in \{1,\ldots,n\}$.  The
\emph{Laplacian} matrix is $ \Lap = \Dout -\Adj$. Note that $
\Lap\ones_n=0 $.  If $ \Bgraph $ is strongly connected, then zero is a
simple eigenvalue of $\Lap$.  $\Bgraph$ is undirected if $ \Lap=\Lap^T
$ and \emph{weight-balanced} if $ \wdout(v) =\wdin(v) $, for all $ v
\in \vertices$. Equivalently, $\Bgraph $ is weight-balanced if and
only if $ \ones_n^T \Lap =0 $ if and only if $ \Lap+\Lap^T $ is
positive semidefinite. Furthermore, if $\Bgraph $ is weight-balanced
and strongly connected, then zero is a simple eigenvalue of $
\Lap+\Lap^T $.  Note that any undirected graph is weight-balanced.
 
\subsection{Zero-sum games}

We recall basic game-theoretic notions following~\cite{TB-GJO:99}.  An
$ n $-player game is a triplet $ \mathbf{G}=(P, \setX,U) $, where $ P
$ is the set of players with $ |P| = n\in \integers_{\geq 2} $, $
\setX=\setX_1\times \ldots \times \setX_n $, $ \setX_i \subset
\real^{d_i} $ is the set of (pure) strategies of player $ v_i \in P $,
$ d_i \in \integerspositive $, and $ U=(u_1,\ldots, u_n) $, where $
\map{u_i}{\setX}{\real} $ is the payoff function of player $ v_i $, $
i\in \{1,\ldots, n\} $.  The game $ \mathbf{G} $ is called a
\emph{zero-sum game} if $ \sum_{i=1}^nu_i = 0 $.  An outcome $ x^* \in
\setX $ is a (pure) \emph{Nash equilibrium} of $ \mathbf{G} $ if for
all $ i \in \{1,\ldots, n\} $ and all $ x_i \in \setX_i $,
\begin{align*}
  u_i(x_i^*, x_{-i}^*)\geq u_i(x_i, x_{-i}^*),
\end{align*}
where $ x_{-i} $ denotes the actions of all players other than~$v_i$.
In this paper, we focus on a class of two-player zero-sum games which
have at least one pure Nash equilibrium as the next result states.

\begin{theorem}\longthmtitle{Minmax theorem}\label{theorem:minmax}
  Let $ \setX_1\subset \real^{d_1} $ and $ \setX_2 \subset
  \real^{d_2}$, $ d_1,d_2\in \integerspositive $, be nonempty,
  compact, and convex. If $ u:\setX_1\times \setX_2 \rightarrow
  \mathbb{R} $ is continuous and the sets $ \{x'\in \setX_1 \ | \
  u(x',y) \geq \alpha\}$ and $ \{y'\in \setX_2 \ | \ u(x,y') \leq
  \alpha\}$ are convex for all $ x \in \setX_1 $, $ y \in \setX_2 $,
  and $ \alpha \in \real $, then
  \[
  \max_x\min_y u(x,y)=\min_y\max_x u(x,y).
  \]
\end{theorem}
Theorem~\ref{theorem:minmax} implies that the game $
\mathbf{G}=(\{v_1,v_2\}, \setX_1\times \setX_2, (u,-u)) $ has a pure
Nash equilibrium.

\myclearpage

\section{Problem statement}\label{sec:problem-statement}

Consider two networks $\network_1$ and $\network_2$ composed of agents
$\{v_1,\dots,v_{n_1}\}$ and agents $\{w_1,\dots,w_{n_2}\}$,
respectively. Throughout this paper, $\network_1 $ and $ \network_2 $
are either connected undirected graphs,
c.f. Section~\ref{sec:nash-undirected}, or strongly connected
weight-balanced digraphs, c.f. Section~\ref{sec:nash-directed}.  Since
the latter case includes the first one, throughout this section, we
assume the latter.  The state of $\network_1$, denoted $x_1$, belongs
to $ \setX_1 \subset \real^{d_1}$, $ d_1 \in \integerspositive
$. Likewise, the state of $\network_2$, denoted $x_2$, belongs to
$\setX_2 \subset \real^{d_2}$, $ d_2 \in \integerspositive $.  In this
paper, we do not get into the details of what these states represent
(as a particular case, the network state could correspond to the
collection of the states of agents in it). In addition, each agent
$v_i $ in $\network_1$ has an estimate $x_1^i \in \real^{d_1}$ of what
the network state is, which may differ from the actual value $x_1$.
Similarly, each agent $w_j $ in $\network_2$ has an estimate $x_2^j
\in \real^{d_2}$ of what the network state is.  Within each network,
neighboring agents can share their estimates.  Networks can also
obtain information about each other. This is modeled by means of a
bipartite directed graph $\networkeng$, called \emph{engagement}
graph, with disjoint vertex sets $\{v_1,\ldots, v_{n_1}\}$ and
$\{w_1,\ldots, w_{n_2}\}$, where every agent has at least one
out-neighbor.  According to this model, an agent in~$\network_1$
obtains information from its out-neighbors in $\networkeng$ about
their estimates of the state of $\network_2$, and vice versa.

\begin{figure}[htb!]
  \centering
  \includegraphics[width=.625\linewidth]{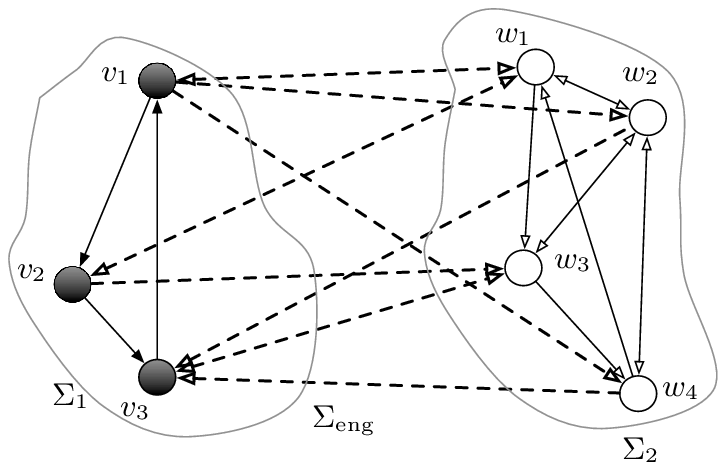}
  \caption{Networks $ \network_1 $ and $ \network_2 $ engaged in a
    strategic scenario.  Both networks are strongly connected and
    weight-balanced, with weights of $ 1 $ on each edge.  Edges which
    correspond to $ \networkeng $ are dashed.}\label{fig:networks}
\end{figure}

Figure~\ref{fig:networks} illustrates this concept.

For each $i \in \{1,\ldots, n_1\}$, let $\map{f_1^i}{\setX_1\times
  \setX_2}{\real} $ be a locally Lipschitz concave-convex function
only available to agent $ v_i \in \network_1 $.  Similarly, let
$\map{f_2^j}{\setX_1\times \setX_2}{\real} $ be a locally Lipschitz
concave-convex function only available to agent $ w_j \in \network_2
$, $ j \in \{1,\ldots,n_2\}$.  The networks $\network_1$ and
$\network_2$ are engaged in a zero-sum game with payoff function
$\map{U}{\setX_1\times \setX_2}{\real} $
\begin{equation}\label{eq:payoff-is-sum}
  U(x_1,x_2)=\sum_{i=1}^{n_1}f_1^{i}(x_1,x_2) =
  \sum_{j=1}^{n_2}f_2^{j}(x_1,x_2), 
\end{equation}
where $\network_1$ wishes to maximize $U$, while $\network_2$ wishes to
minimize it. The objective of the networks is therefore to settle upon
a Nash equilibrium, i.e., to solve the following maxmin problem
\begin{equation}\label{eq:problem_1}
  \max_{x_1\in \setX_1} \min_{x_2 \in \setX_2} U(x_1,x_2) .
\end{equation}
We refer to the this zero-sum game as the \emph{$2$-network zero-sum
  game} and denote it by $\gameforiginal$.  We assume that $ \setX_1
\subset \real^{d_1}$ and $ \setX_2 \subset \real^{d_2} $ are compact
convex. For convenience, let $ \xnet_1=(x_1^1, \ldots, x_1^{n_1})^T $
and $ \xnet_2=(x_2^1, \ldots, x_2^{n_2})^T $ denote vector of agent
estimates about the state of the respective networks.

\begin{remark}\longthmtitle{Power allocation in communication channels
    in the presence of adversaries}\label{remark:application}
  {\rm Here we present an example from communications inspired
    by~\citep[Section 5.5.3]{SB-LV:04}.  Consider $n$ Gaussian
    communication channels, each with signal power $p_i \in
    \realnonnegative$ and noise power $\eta_i \in \realnonnegative$,
    for $i \in \until{n}$.  The capacity of each channel is
    proportional to $ \log(1+\beta p_i/(\sigma_i+\eta_i))$, where $
    \beta \in \realpositive $ and $ \sigma_i >0 $ is the receiver
    noise.  Note that capacity is concave in $p_i$ and convex in
    $\eta_i$.  Both signal and noise powers must satisfy a budget
    constraint, i.e., $ \sum_{i=1}^{n} p_i=P $ and $
    \sum_{i=1}^{n}\eta_i = C$, for some given $ P, C \in
    \realpositive$. Two networks of $n$ agents are involved in this
    scenario, one, $\Sigma_1$, selecting signal powers to maximize
    capacity, the other one, $\Sigma_2$, selecting noise powers to
    minimize it. The network $\Sigma_1$ has decided that $m_1$
    channels will have signal power $x_1$, while $n-1-m_1$ will have
    signal power $x_2$. The remaining $n$th channel has its power
    determined to satisfy the budget constraint, i.e., $P-m_1 x_1 -
    (n-1-m_1)x_2$. Likewise, the network $\Sigma_2$ does something
    similar with $m_2$ channels with noise power $y_1$, $n-1-m_2$
    channels with noise power $y_2$, and one last channel with noise
    power $C-m_2 y_1 - (n-1-m_2) y_2$. Each network is aware of the
    partition made by the other one.  The individual objective
    function of the two agents (one from~$\Sigma_1$, the other from
    $\Sigma_2$) making decisions on the power levels of the $i$th
    channel is the channel capacity itself. For $i \in \until{n-1}$,
    this takes the form
    \begin{align*}
      f^i (x,y)= \log \Big(1+ \frac{\beta x_a}{\sigma_i+y_b} \Big),
    \end{align*}
    for some $a,b \in \{1,2\}$. Here $x=(x_1,x_2)$ and
    $y=(y_1,y_2)$. For $i=n$, it takes instead the form
    \begin{align*}
      f^{n}(x,y) = \log \Big(1+ \frac{\beta (P-m_1 x_1 -
        (n-1-m_1)x_2)}{\sigma_{n}+C-m_2 y_1 - (n-1-m_2) y_2} \Big).
    \end{align*}
    Note that $ \sum_{i=1}^{n}f^i(x,y) $ is the total capacity of the
    $n$ communication channels.  } \oprocend
\end{remark}


\subsection{Reformulation of the $2$-network zero-sum game}\label{sec:reformulation}

In this section, we describe how agents in each network use the
information obtained from their neighbors to compute the value of
their own objective functions. Based on these estimates, we introduce a
reformulation of the $ \gameforiginal $ which is instrumental for
establishing some of our results.

Each agent in $\network_1$ has a locally Lipschitz, concave-convex
function $ \map{\fLift_1^i}{\real^{d_1}\times \real^{d_2n_2}}{\real} $
with the properties:
\begin{description}
\item[(Extension of own payoff function):] for any $ x_1 \in
  \real^{d_1}$, $x_2 \in \real^{d_2}$,
  \begin{subequations}\label{eq:properties}
    \begin{equation}\label{eq:property_1}
      \fLift_1^i(x_1,\ones_{n_2}\otimes x_2)=f_1^i(x_1,x_2).
    \end{equation}
    
  \item[(Distributed over $\networkeng$):] there exists $
    \map{\fdist_1^i}{\real^{d_1} \times
      \real^{d_2|\Nin_{\networkeng}(v_i)|}}{\real}$ such that, for any
    $ x_1 \in \real^{d_1}$ $ \xnet_2 \in \real^{d_2n_2} $,
    \begin{equation}\label{eq:property_2}
      \fLift_1^i(x_1,\xnet_2) = \fdist_1^i(x_1,\pi_1^i(\xnet_2)) ,
    \end{equation}
  \end{subequations}
  with $
  \map{\pi_1^i}{\real^{d_2n_2}}{\real^{d_2|\Nout_{\networkeng}(v_i)|}}$
  the projection of $ \xnet_2 $ to the values received by $ v_i $ from
  its out-neighbors in~$ \networkeng $.
\end{description}
Equation~\eqref{eq:property_1} states the fact that, when the
estimates of all neighbors of an agent in the opponent's network
agree, its evaluation should coincide with this
estimate. Equation~\eqref{eq:property_2} states the fact that agents
can only use the information received from their neighbors in the
interaction topology to compute their new estimates.

Each agent in $\network_2$ has a function $
\map{\fLift_2^j}{\real^{d_1n_1}\times \real^{d_2}}{\real} $ with
similar properties.  The collective payoff functions of the two
networks are
\begin{subequations}\label{eq:payoff-is-sum_2}
  \begin{align}
    \ULift_1(\xnet_1,\xnet_2) &=
    \sum_{i=1}^{n_1}\fLift_1^i(x_1^i,\xnet_2),
    \\
    \ULift_2(\xnet_1,\xnet_2) &=
    \sum_{j=1}^{n_2}\fLift_2^j(\xnet_1,x_2^j).
  \end{align}  
\end{subequations}
In general, the functions $ \ULift_1 $ and $ \ULift_2 $ need not be
the same.  However, $ \ULift_1 (\ones_{n_1} \otimes x_1,\ones_{n_1}
\otimes x_2) = \ULift_2 (\ones_{n_1} \otimes x_1,\ones_{n_1} \otimes
x_2)$, for any $x_1 \in \real^{d_1}$, $x_2 \in \real^{d_2}$.  When
both functions coincide, the next result shows that the original game
can be lifted to a (constrained) zero-sum game.

\begin{lemma}\longthmtitle{Reformulation of the $2$-network zero-sum
    game}\label{le:equiv-Utilde}
  Assume that the individual payoff functions
  $\{\fLift^i_1\}_{i=1}^{n_1}$, $\{\fLift^j_2\}_{j=1}^{n_2}$
  satisfying~\eqref{eq:properties} are such that the network payoff
  functions defined in~\eqref{eq:payoff-is-sum_2} satisfy $ \ULift_1 =
  \ULift_2 $, and let~$ \ULift $ denote this common function.  Then,
  the problem~\eqref{eq:problem_1} on $\real^{d_1}\times\real^{d_2} $
  is equivalent to the following problem on $\real^{n_1d_1}\times
  \real^{n_2d_2}$,
  \begin{align}\label{eq:problem_stat_2}
    & \max_{\xnet_1\in \setX_1^{n_1}} \min_{\xnet_2 \in \setX_2^{n_2}}
    \ULift(\xnet_1,\xnet_2), \nonumber \\
    & \qquad \text{subject to} \quad \Blap_1 \xnet_1 =
    \zeros_{n_1d_1}, \quad \Blap_2 \xnet_2 = \zeros_{n_2d_2} ,
  \end{align}
  with $ \Blap_\ell= \Lap_\ell \otimes \identity{d_\ell}$ and $
  \Lap_\ell $ the Laplacian of $\network_\ell$,~$ \ell \in \{1, 2\} $.
\end{lemma}
\begin{pf*}{Proof.}
  The proof follows by noting that (i) $\ULift (\ones_{n_1}\otimes
  x_1, \ones_{n_2}\otimes x_2) = U(x_1,x_2)$ for all $ x_1 \in
  \real^{d_1} $ and $ x_2 \in \real^{d_2} $ and (ii) since $\Bgraph_1$
  and $\Bgraph_2$ are strongly connected, $\Blap_1 \xnet_1 =
  \zeros_{n_1d_1} $ and $\Blap_2\xnet_2 = \zeros_{n_2d_2} $ iff
  $\xnet_1 = \ones_{n_1} \otimes x_1 $ and $\xnet_2 = \ones_{n_2}
  \otimes x_2 $ for some $x_1 \in \real^{d_1}$ and $x_2 \in
  \real^{d_2}$. \qed
\end{pf*}

\begin{remark}\longthmtitle{Restrictions on extensions}
  {\rm The assumption of Lemma~\ref{le:equiv-Utilde} does not hold in
    general for all sets of extensions
    satisfying~\eqref{eq:property_1} and~\eqref{eq:property_2}.  If
    the interaction topology is one-to-one (i.e., both networks have
    the same number of agents, the interaction topology is undirected,
    and each agent in the first network obtains information only from
    one agent in the opposing network), the natural extensions satisfy
    the assumption.  Example~\ref{ex:example} later provides yet
    another instance of a different nature.  In general, determining
    if it is always possible to choose the extensions in such a way
    that the assumption holds is an open problem. }  \oprocend
\end{remark}

We denote by $\gamef$ the constrained zero-sum game defined
by~\eqref{eq:problem_stat_2} and refer to this situation by saying
that $\gameworiginal$ can be lifted to~$\gamew$.  Our objective is to
design a coordination algorithm that is implementable with the
information that agents in $\network_1$ and $\network_2$ possess and
leads them to find a Nash equilibrium of $ \gamew $, which corresponds
to a Nash equilibrium of $ \gameworiginal $ by
Lemma~\ref{le:equiv-Utilde}.  Achieving this goal, however, is
nontrivial because individual agents, not networks themselves, are the
decision makers. From the point of view of agents in each network, the
objective is to agree on the states of both their own network and the
other network, and that the resulting states correspond to a Nash
equilibrium of $ \gameworiginal $.

%

The function $\ULift$ is locally Lipschitz and
concave-convex. Moreover, from Lemma~\ref{le:finite_sum_gen}, the
elements of $\partialvector_{\xnet_1} \ULift(\xnet_1,\xnet_2)$ are of
the form
\begin{align*}
  \gLift_{(\xnet_1,\xnet_2)} = (\gLift^1_{(x^1_1, \xnet_2)}, \ldots,
  \gLift^n_{(x^n_1, \xnet_2)}) \in \partialvector_{\xnet_1}
  \ULift(\xnet_1,\xnet_2),
\end{align*}
where $ \gLift^i_{(x^i_1, \xnet_2)} \in \partial_{x_1}
\fLift^i_1(x^i_1, \xnet_2) $, for $ i \in \{1,\ldots, n_1\}$.  Note
that, because of~\eqref{eq:property_2}, we have $\partial_{x_1}
\fLift^i_1(x^i_1, \ones_{n_2}\otimes x_2 ) = \partial_{x_1}
f^i_1(x^i_1, x_2) $.  A similar reasoning can be followed to describe
the elements of $ \partialvector_{\xnet_2} \ULift(\xnet_1,\xnet_2) $.
Next, we present a characterization of the Nash equilibria of $
\gamew$, instrumental for proving some of our upcoming results.

\begin{proposition}\longthmtitle{Characterization of the Nash
    equilibria of $\gamew $}\label{prop:equiv-F}
  For $\network_1$, $\network_2$ strongly connected and
  weight-balanced, define $ F_1 $ and $ F_2 $ by
  \begin{align*}
    F_1(\xnet_1,\znet_1,\xnet_2)&=-\ULift(\xnet_1,\xnet_2) +
    \xnet_1^T\Blap_1 \znet_1 + \frac{1}{2}\xnet_1^T\Blap_1\xnet_1,
    \\
    F_2(\xnet_2,\znet_2,\xnet_1)& =\ULift(\xnet_1,\xnet_2) +
    \xnet_2^T\Blap_2 \znet_2 + \frac{1}{2}\xnet_2^T\Blap_2\xnet_2.
  \end{align*}
  Then, $F_1$ and $F_2$ are convex in their first argument, linear in
  their second one, and concave in their third one.  Moreover, assume
  $(\xnet^*_1,\znet^*_1,\xnet^*_2,\znet^*_2) $ satisfies the following
  \emph{saddle property} for $(F_1,F_2)$: $(\xnet^*_1,\znet^*_1) $ is
  a saddle point of $(\xnet_1,\znet_1) \mapsto F_1(\xnet_1,\znet_1,
  \xnet^*_2) $ and $(\xnet^*_2,\znet^*_2) $ is a saddle point of
  $(\xnet_2,\znet_2) \mapsto F_2(\xnet_2,\znet_2, \xnet^*_1) $. Then,
  \begin{enumerate}
  \item $(\xnet^*_1,\znet^*_1+\ones_{n_1}\otimes
    a_1,\xnet^*_2,\znet^*_2+\ones_{n_2}\otimes a_2) $ satisfies the
    saddle property for $(F_1,F_2) $ for any $ a_1\in \real^{d_1} $, $
    a_2 \in\real^{d_2} $, and
  \item $(\xnet^*_1,\xnet^*_2) $ is a Nash equilibrium of $ \gamew $.
  \end{enumerate}
  Furthermore,
  \begin{enumerate}
  \item[(iii)] if $ (\xnet^*_1,\xnet^*_2) $ is a Nash equilibrium of $
    \gamew $ then there exists $ \znet^*_1,\znet^*_2 $ such that
    $(\xnet^*_1,\znet^*_1,\xnet^*_2,\znet^*_2) $ satisfies the saddle
    property for $(F_1,F_2) $.
  \end{enumerate}
\end{proposition}

\begin{pf*}{Proof.}
  The statement~(i) is immediate. For~(ii), since
  $(\xnet^*_1,\znet^*_1,\xnet^*_2,\znet^*_2) $ satisfies the saddle
  property, and the networks are strongly connected and
  weight-balanced, we have $\xnet^*_1 = \ones_{n_1} \otimes x^*_1 $, $
  x^*_1 \in \real^{d_1} $, $\xnet^*_2 = \ones_{n_2}\otimes x^*_2 $, $
  x^*_2 \in \real^{d_2} $, $\Blap_1 \znet^*_1 \in
  -\partialvector_{\xnet_1} \ULift(\xnet^*_1,\xnet^*_2) $, and
  $\Blap_2 \znet^*_2 \in
  \partialvector_{\xnet_2} \ULift (\xnet^*_1,\xnet^*_2) $. Thus there
  exist $ g^i_{1,(x^*_1,x^*_2)} \in \partial_{x_1} f_1^i(x^*_1,x^*_2)
  $, $ i \in \{1,\ldots, n_1\} $, and $ g^j_{2,(x^*_1,x^*_2)} \in
  -\partial_{x_2} f_2^j(x^*_1,x^*_2)$, $ j \in \{1,\ldots, n_2\} $,
  such that
  \begin{align*}
    \Blap_1\znet_1^*& =(g^1_{1,(x^*_1,x^*_2)}, , \ldots,
    g^n_{1,(x^*_1,x^*_2)}))^T, \ \mathrm{and} \\
    \Blap_2 \znet_2^*&=(g^1_{2,(x^*_1,x^*_2)}, , \ldots,
    g^n_{2,(x^*_1,x^*_2)}))^T.
  \end{align*}
  Noting that, for $ \ell \in \{1,2\} $, $ (\ones_{n_\ell}^T\otimes
  \identity{d_\ell}) \Blap_\ell = (\ones_{n_\ell}^T\otimes \identity{d_\ell})
  (\Lap_\ell \otimes \identity{d_\ell}) = \ones_{n_\ell}^T \Lap \otimes
  \identity{d_\ell} = \zeros_{d_\ell \times d_\ell n_\ell}$, we deduce that $
  \sum_{i=1}^{n_1} g^i_{1,(x^*_1,x^*_2)}=\mathbf{0}_{d_1} $ and $
  \sum_{j=1}^{n_2} g^j_{2,(x^*_1,x^*_2)}= \mathbf{0}_{d_2} $, i.e., $
  (\xnet^*_1,\xnet^*_2) $ is a Nash equilibrium.
  Finally for proving~(iii), note that $\xnet_1^* = \ones_{n_1}\otimes
  x_1^*$ and $\xnet_2^* = \ones_{n_2}\otimes x_2^*$. The result
  follows then from the fact that $0 \in \partial_{x_1}U(x_1^*,x_2^*)$
  and $0 \in \partial_{x_2}U(x_1^*,x_2^*)$ implies that there exists
  $\znet_1^*\in \real^{n_1d_1} $ and $\znet_2^*\in \real^{n_2d_2} $
  with $\Blap_1 \znet_1^*
  \in \partial_{\xnet_1}\ULift(\xnet_1^*,\xnet_2^*)$ and $\Blap_2
  \znet_2^* \in -\partial_{\xnet_2}\ULift(\xnet_1^*,\xnet_2^*)$.
  \qed
\end{pf*}

\myclearpage
\section{Distributed convergence to Nash equilibria for undirected
  topologies}\label{sec:nash-undirected}

In this section, we introduce a distributed dynamics which
solves~\eqref{eq:problem_stat_2} when $ \network_1 $ and $ \network_2$
are undirected. In particular, we design gradient dynamics to find
points with the saddle property for $(F_1,F_2) $ prescribed by
Proposition~\ref{prop:equiv-F}. Consider the set-valued dynamics
$\setmap{\undirectedNash}{(\real^{d_1n_1})^2 \times
  (\real^{d_2n_2})^2} {(\real^{d_1n_1})^2\times (\real^{d_2n_2})^2} $,
\begin{subequations}\label{eq:Nash_dynamics-undirected}
  \begin{align}
    \dot \xnet_1 + \Blap_1 \xnet_1 + \Blap_1 \znet_1&
    \in \partialvector_{\xnet_1} \ULift
    (\xnet_1,\xnet_2), \label{eq:Nash_dynamics-undirected-a}
    \\
    \dot \znet_1 &= \Blap_1
    \xnet_1, \label{eq:Nash_dynamics-undirected-b}
    \\
    \dot \xnet_2 + \Blap_2\xnet_2+\Blap_2\znet_2 & \in
    - \partialvector_{\xnet_2}\ULift
    (\xnet_1,\xnet_2), \label{eq:Nash_dynamics-undirected-c}
    \\
    \dot \znet_2 &= \Blap_2
    \xnet_2, \label{eq:Nash_dynamics-undirected-d}
  \end{align}
\end{subequations}
where $ \xnet_\ell,\znet_\ell \in \real^{n_\ell d_\ell} $, $ \ell \in
\{1,2\} $.
Note
that~\eqref{eq:Nash_dynamics-undirected-a}-\eqref{eq:Nash_dynamics-undirected-b}
and~\eqref{eq:Nash_dynamics-undirected-c}-\eqref{eq:Nash_dynamics-undirected-d}
correspond to saddle-point dynamics of $F_1$ in $(\xnet_1,\znet_1)$
and $F_2$ in $(\xnet_2,\znet_2)$, respectively.  Local solutions to
this dynamics exist by virtue of Lemmas~\ref{le:gradient_properties}
and~\ref{le:solution}.  We characterize next its asymptotic
convergence properties.

\begin{theorem}\longthmtitle{Distributed  convergence to Nash equilibria
    for undirected networks}\label{theorem:Nash-undirected}
  Consider the zero-sum game $ \gameforiginal $, with $ \network_1 $
  and $ \network_2 $ connected undirected graphs, $ \setX_1 \subset
  \real^{d_1} $, $ \setX_2 \subset \real^{d_2} $ compact and convex,
  and $\map{U}{\setX_1 \times \setX_2}{\real} $ strictly
  concave-convex and locally Lipschitz. Assume $\gameworiginal$ can be
  lifted to $\gamew$.  Then, the projection onto the first and third
  components of the solutions of~\eqref{eq:Nash_dynamics-undirected}
  asymptotically converge to agreement on the Nash equilibrium of
  $\gameworiginal$.
\end{theorem}

\begin{pf*}{Proof.}
  Throughout this proof, since property~\eqref{eq:property_2} holds,
  without loss of generality and for simplicity of notation, we assume
  that agents in $ \network_1 $ have access to $ \xnet_2 $ and,
  similarly, agents in $ \network_2 $ have access to $ \xnet_1 $.  By
  Theorem~\ref{theorem:minmax}, a solution
  to~\eqref{eq:payoff-is-sum_2} exists. By the strict
  concavity-convexity properties, this solution is, in fact, unique.
  Let us denote this solution by $ \xnet_1^*=\ones_{n_1} \otimes
  x^*_1$ and $ \xnet_2^*=\ones_{n_2} \otimes x^*_2 $. By
  Proposition~\ref{prop:equiv-F}(iii), there exists $\znet_1^*$ and $
  \znet_2^* $ such that $ (\xnet_1^*,\znet_1^*, \xnet_2^*,\znet_2^*)
  \in \equilibria{\undirectedNash} $.  First, note that given any
  initial condition $ (\xnet^0_1,\znet^0_1, \xnet^0_2,\znet^0_2) \in
  (\real^{n_1d_1})^2 \times (\real^{n_2d_2})^2 $, the set
  \begin{align}\label{eq:S_z0}
    W_{\znet^0_1,\znet^0_2} & = \{(\xnet_1,\znet_1,\xnet_2,\znet_2) \
    | \nonumber
    \\
    & \ (\ones_{n_\ell}^T \otimes \identity{d_\ell}) \znet_\ell =
    (\ones_{n_\ell}^T \otimes \identity{d_\ell}) \znet^0_\ell, \; \ell
    \in \{1,2\} \}
  \end{align}
  is strongly positively invariant
  under~\eqref{eq:Nash_dynamics-undirected}. Consider the function
  $\map{V}{(\real^{d_1n_1})^2 \times
    (\real^{d_2n_2})^2}{\realnonnegative}$ defined by
  \begin{align*}
    V(&\xnet_1,\znet_1,\xnet_2,\znet_2) \\
    &= \frac{1}{2}(\xnet_1-\xnet_1^*)^T(\xnet_1-\xnet_1^*) +
    \frac{1}{2}(\znet_1-\znet_1^*)^T(\znet_1-\znet_1^*)
    \\
    & \quad + \frac{1}{2}(\xnet_2-\xnet_2^*)^T(\xnet_2-\xnet_2^*) +
    \frac{1}{2}(\znet_2-\znet_2^*)^T(\znet_2-\znet_2^*).
  \end{align*}
  The function $ V $ is smooth. Next, we examine its set-valued Lie
  derivative along $ \undirectedNash $. Let $ \xi \in
  \SetLie_{\undirectedNash}V(\xnet_1,\znet_1,\xnet_2,\znet_2) $.  By
  definition, there exists $ v \in
  \undirectedNash(\xnet_1,\znet_1,\xnet_2,\znet_2) $, given by
  \begin{align*}
    v=(&-\Blap_1 \xnet_1 -\Blap_1 \znet_1 +g_{1,(\xnet_1,\xnet_2)}, \\
    & -\Blap_2 \xnet_2 -\Blap_2 \znet_2 -g_{2,(\xnet_1,\xnet_2)} ,
    \Blap_1\xnet_1, \Blap_2\xnet_2),
  \end{align*}
  where $ g_{1,(\xnet_1,\xnet_2)} \in \partialvector_{\xnet_1}
  U(\xnet_1,\xnet_2) $ and $ g_{2,(\xnet_1,\xnet_2)}
  \in \partialvector_{\xnet_2}U (\xnet_1,\xnet_2) $, such that
  \begin{align*}
    \xi & =v\cdot \nabla V(\xnet_1,\znet_1,\xnet_2,\znet_2)
    \\
    &= (\xnet_1-\xnet_1^*)^T(-\Blap_1 \xnet_1 -\Blap_1 \znet_1
    +g_{1,(\xnet_1,\xnet_2)}) \\
    &+ (\xnet_2-\xnet_2^*)^T(-\Blap_2 \xnet_2
    -\Blap_2 \znet_2 -g_{2,(\xnet_1,\xnet_2)} )
    \\
    & \quad + (\znet_1-\znet_1^*)^T \Blap_1\xnet_1 +
    (\znet_2-\znet_2^*)^T\Blap_2\xnet_2.
  \end{align*}  
  Note that $ -\Blap_1 \xnet_1 -\Blap_1 \znet_1
  +g_{1,(\xnet_1,\xnet_2)}\in -\partialvector_{\xnet_1}
  F_1(\xnet_1,\znet_1,\xnet_2)$,
  $\Blap_1\xnet_1\in \partialvector_{\znet_1}
  F_1(\xnet_1,\znet_1,\xnet_2) $, $ -\Blap_2 \xnet_2 -\Blap_2 \znet_2
  -g_{2,(\xnet_1,\xnet_2)} \in -\partialvector_{\xnet_2}
  F_2(\xnet_1,\znet_2,\xnet_2) $, and $
  \Blap_2\xnet_2\in \partialvector_{\znet_2}F_2(\xnet_2,\znet_2,\xnet_1)$.
  Using the first-order convexity property of $ F_1 $ and $ F_2 $ in
  their first two arguments, one gets
  \begin{align*}
    \xi & \leq
    F_1(\xnet_1^*,\znet_1,\xnet_2)-F_1(\xnet_1,\znet_1,\xnet_2)
    +F_2(\xnet_2^*,\znet_2,\xnet_1)\\
    & -F_2(\xnet_2,\znet_2,\xnet_1)
    +F_1(\xnet_1,\znet_1,\xnet_2)-F_1(\xnet_1,\znet^*_1,\xnet_2)\\    
    &+F_2(\xnet_2,\znet_2,\xnet_1)-F_2(\xnet_2,\znet_2^*,\xnet_1).
  \end{align*}
  Expanding each term and using the fact that $
  (\xnet_1^*,\znet_1^*,\xnet_2^*,\znet_2^*) $ $ \in
  \equilibria{\undirectedNash}$, we simplify this inequality~as
  \begin{align*}
    \xi & \leq -\ULift(\xnet_1^*,\xnet_2) +\ULift(\xnet_1,\xnet^*_2)
    -\znet^*_1\Blap_1\xnet_1
    \\
    & \quad -\frac{1}{2}\xnet_1\Blap_1\xnet_1
    -\znet^*_2\Blap_2\xnet_2-\frac{1}{2}\xnet_2\Blap_2\xnet_2 .
  \end{align*}
  By rearranging, we thus have
  \begin{align*}
    \xi \leq
    -F_2(\xnet_2,\znet_2^*,
    \xnet_1^*)-F_1(\xnet_1,\znet_1^*,\xnet_2^*) .
  \end{align*}
  Next, since $
  F_2(\xnet_1^*,\znet_2^*,\xnet_2^*)+F_1(\xnet_2^*,\znet_2^*,\xnet_1^*)
  =0$, we have
  \begin{align*}
    \xi & \leq F_1(\xnet_1^*,\znet_1^*,\xnet_2^*)
    -F_1(\xnet_1,\znet_1^*,\xnet_2^*)\\
    & \quad +
    F_2(\xnet_2^*,\znet_2^*,\xnet_1^*)-F_2(\xnet_2,\znet_2^*,\xnet^*_1),
  \end{align*}
  yielding that $ \xi \leq 0 $. As a result, 
  \[ \max
  \SetLie_{\undirectedNash}V (\xnet_1,\znet_1, \xnet_2,\znet_2) \leq 0.
  \]
  As a by-product, we conclude that the trajectories
  of~\eqref{eq:Nash_dynamics-undirected} are bounded.
  By virtue of the set-valued version of the LaSalle Invariance
  Principle, cf. Theorem~\ref{th:laSalle}, any trajectory
  of~\eqref{eq:Nash_dynamics-undirected} starting from an initial
  condition $ (\xnet_1^0,\znet_1^0,\xnet_2^0,\znet_2^0) $ converges to
  the largest positively invariant set $ M $ in $
  S_{\undirectedNash,\Lyap} \cap \Lyap^{-1}(\leq
  V(\xnet_1^0,\znet_1^0,\xnet_2^0,\znet_2^0))$.  Let
  $(\xnet_1,\znet_1,\xnet_2,\znet_2) \in M$. Because $M \subset
  S_{\undirectedNash,\Lyap} $, then $
  F_1(\xnet_1^*,\znet_1^*,\xnet_2^*)
  -F_1(\xnet_1,\znet_1^*,\xnet_2^*)=0, $~i.e.,
  \begin{equation}\label{eq:F_1-ineq}
    -\ULift(\xnet_1^*,\xnet_2^*) +
    \ULift(\xnet_1,\xnet_2^*)-\xnet_1^T\Blap_1\znet_1^* -
    \frac{1}{2}\xnet_1^T\Blap_1\xnet_1 = 0.
  \end{equation}
  Define now $ \map{G_1}{\real^{n_1d_1} \times \real^{n_1d_1} \times
    \real^{n_2d_2}}{\real} $ by $G_1(\xnet_1,\znet_1,\xnet_2) =
  F_1(\xnet_1,\znet_1,\xnet_2)-\frac{1}{2}\xnet_1^T\Blap_1\xnet_1$.
  $G_1$ is convex in its first argument and linear in its
  second. Furthermore, for fixed $\xnet_2$, the map $(\xnet_1,\znet_1)
  \mapsto G_1(\xnet_1,\znet_1,\xnet_2)$ has the same saddle points as
  $(\xnet_1,\znet_1) \mapsto F_1(\xnet_1,\znet_1,\xnet_2)$.  As a
  result, $ G_1(\xnet_1^*,\znet_1^*,\xnet_2^*) -
  G_1(\xnet_1,\znet^*_1,\xnet_2^*) \leq 0 $, or equivalently, $
  -\ULift(\xnet_1^*,\xnet_2^*) + \ULift(\xnet_1,\xnet_2^*) - \xnet_1^T
  \Blap_1 \znet_1^* \leq 0 $.  Combining this
  with~\eqref{eq:F_1-ineq}, we have that $ \Blap_1\xnet_1=0 $ and $
  -\ULift(\xnet_1^*,\xnet_2^*)+\ULift(\xnet_1,\xnet_2^*)= 0 $.  Since
  $ \ULift $ is strictly concave in its first argument $
  \xnet_1=\xnet_1^* $. A similar argument establishes that $
  \xnet_2=\xnet_2^* $.  Using now the fact that $ M $ is weakly
  positively invariant, one can deduce that $ \Blap_\ell \znet_\ell
  \in -\partial_{\xnet_\ell} \ULift(\xnet_1,\xnet_2) $, for $ \ell
  \in\{1,2\} $, and thus $ (\xnet_1,\znet_1, \xnet_2,\znet_2)\in
  \equilibria{\undirectedNash} $. \qed
\end{pf*}

\myclearpage
\section{Distributed convergence to Nash equilibria for directed
  topologies}\label{sec:nash-directed}
Interestingly, the saddle-point
dynamics~\eqref{eq:Nash_dynamics-undirected} fails to converge when
transcribed to the directed network setting.  This observation is a
consequence of the following result, which studies the stability of
the linearization of the dynamics~\eqref{eq:Nash_dynamics-undirected},
when the payoff functions have no contribution to the linear part.

\begin{lemma}\longthmtitle{Necessary condition for the convergence
    of~\eqref{eq:Nash_dynamics-undirected} on
    digraphs}\label{lemma:necessary_linear_directed}
  Let $\network_\ell $ be strongly connected and $ f^i_\ell=0$, $ i
  \in \{1,\ldots, n_\ell\}$, for $ \ell \in\{1,2\} $. Then, the set of
  network agreement configurations $\SS= \setdef{(\ones_{n_1} \otimes
    x_1, \ones_{n_1} \otimes z_1,\ones_{n_2} \otimes x_2, \ones_{n_2}
    \otimes z_2) \in (\real^{n_1 d_1})^2 \times (\real^{n_2
      d_2})^2}{x_\ell,z_\ell \in \real^{d_\ell}, \ell \in \{1,2\} }$,
  is stable under~\eqref{eq:Nash_dynamics-undirected} iff, for any
  nonzero eigenvalue $\lambda$ of the Laplacian $ \Lap_\ell $, $ \ell
  \in \{1,2\} $, one has $\sqrt{3}|\impart (\lambda)| \le \realpart
  (\lambda) $.
\end{lemma}
\begin{pf*}{Proof.}
  In this case,~\eqref{eq:Nash_dynamics-undirected} is linear with
  matrix
  \begin{align}
    \label{eq:auxx}
    \begin{pmatrix}
      \left(\begin{smallmatrix} -1 & -1\\1 & 0\end{smallmatrix}
      \right)\otimes \Blap_1
      & 0 \\
      0 & \left(\begin{smallmatrix} -1 & -1\\1 & 0\end{smallmatrix}
      \right)\otimes \Blap_2
    \end{pmatrix}
  \end{align}
  and has $\SS$ as equilibria. The eigenvalues of~\eqref{eq:auxx} are
  of the form $\lambda_\ell \, \big(\frac{-1}{2}\pm\frac{\sqrt{3}}{2}i
  \big)$, with $\lambda_\ell$ an eigenvalue of $\Blap_\ell$, for $\ell
  \in \{1,2\} $ (since the eigenvalues of a Kronecker product are the
  product of the eigenvalues of the corresponding matrices). Since
  $\Blap_\ell = \Lap_\ell \otimes \identity{d_\ell}$, each eigenvalue
  of $\Blap_\ell$ is an eigenvalue of $ \Lap_\ell $. The result
  follows by noting that $ \realpart \big(\lambda_\ell
  \big(\frac{-1}{2}\pm\frac{\sqrt{3}}{2}i \big)\big) = \frac{1}{2} (
  \mp \sqrt{3}\impart(\lambda_\ell)-\realpart(\lambda_\ell))$.\qed
\end{pf*}

It is not difficult to construct examples of strictly concave-convex
functions that have zero contribution to the linearization
of~\eqref{eq:Nash_dynamics-undirected} around the solution. Therefore,
such systems cannot be convergent if they fail the necessary condition
identified in Lemma~\ref{lemma:necessary_linear_directed}.
The counterexample provided in our recent
paper~\citep{BG-JC:12-tac} of strongly connected, weight-balanced
digraphs that do not meet the stability criterium of
Lemma~\ref{lemma:necessary_linear_directed} is therefore valid in this
context too.

From here on, we assume that the payoff functions are
differentiable. We elaborate on the reasons for this assumption in
Remark~\ref{re:payoff-function-properties} later. Motivated by the
observation made in Lemma~\ref{lemma:necessary_linear_directed}, we
introduce a parameter $ \alpha \in \realpositive $ in the dynamics
of~\eqref{eq:Nash_dynamics-undirected} as
\begin{subequations}\label{eq:Nash_dynamics-directed}
  \begin{align}
    \dot \xnet_1 + \alpha \Blap_1\xnet_1 + \Blap_1\znet_1 & = \gradient
    \ULift(\xnet_1,\xnet_2),
    \\
    \dot \znet_1 & = \Blap_1 \xnet_1,\\
    \dot \xnet_2 + \alpha \Blap_2\xnet_2 + \Blap_2s\znet_2 & = -\gradient
    \ULift(\xnet_1,\xnet_2),
    \\
    \dot \znet_2 & = \Blap_2 \xnet_2.
  \end{align}  
\end{subequations}
We next show that a suitable choice of $\alpha$ makes the dynamics
convergent to the Nash equilibrium.

\begin{theorem}\longthmtitle{Distributed convergence to Nash
    equilibria for directed networks}\label{theorem:Nash-directed}
  Consider the zero-sum game $ \gameforiginal $, with $ \network_1 $
  and $ \network_2 $ strongly connected and weight-balanced digraphs,
  $ \setX_1 \subset \real^{d_1} $, $ \setX_2 \subset \real^{d_2} $
  compact and convex, and $\map{U}{\setX_1 \times \setX_2}{\real} $
  strictly concave-convex and differentiable with globally Lipschitz
  gradient. Assume $\gameworiginal$ can be lifted to $\gamew$ such
  that $\ULift$ is differentiable and has a globally Lipschitz
  gradient.  Define $ \map{h}{\realpositive}{\real} $~by
  \begin{align}\label{eq:h}
    h(r) =& \frac{1}{2}\smallnonzeroeig^{\min} \Big(\sqrt{\Big(
      \frac{r^4+3r^2+2}{r}\Big)^2-4} -\frac{r^4+3r^2+2}{r} \Big)\nonumber \\
    &+\frac{K r^2}{(1+r^2)},
  \end{align}
  where $ \smallnonzeroeig^{\min} =
  \min_{\ell=1,2}\{\smallnonzeroeig(\Lap_\ell+\Lap_\ell^T)\} $, $
  \smallnonzeroeig(\cdot) $ denotes the smallest non-zero eigenvalue
  and $ K \in \realpositive $ is the Lipschitz constant of the
  gradient of $ \ULift $. Then there exists $ \beta^* \in
  \realpositive $ with $ h(\beta^*)=0 $ such that for all $ 0 < \beta
  < \beta^* $, the projection onto the first and third components of
  the solutions of~\eqref{eq:Nash_dynamics-directed} with $
  \alpha=\frac{\beta^2+2}{\beta} $ asymptotically converge to
  agreement on the Nash equilibrium of $\gameworiginal$.
\end{theorem}
  
\begin{pf*}{Proof.}
  Similarly to the proof of Theorem~\ref{theorem:Nash-undirected}, we
  assume, without loss of generality, that agents in $ \network_1 $ have
  access to $ \xnet_2 $ and agents in $ \network_2 $ to $ \xnet_1 $.
  For convenience, we denote the dynamics described
  in~\eqref{eq:Nash_dynamics-directed} by $
  \map{\directedNash}{(\real^{d_1n_1})^2 \times (\real^{d_2n_2})^2}
  {(\real^{d_1n_1})^2 \times (\real^{d_2n_2})^2} $.  Let $
  (\xnet_1^0,\znet_1^0,\xnet_2^0,\znet_2^0) $ be an arbitrary initial
  condition.  Note that the set $ W_{\znet^0_1,\znet^0_2} $ defined
  by~\eqref{eq:S_z0} is invariant under the evolutions
  of~\eqref{eq:Nash_dynamics-directed}.  By an argument similar to the
  one in the proof of Theorem~\ref{theorem:Nash-undirected}, there
  exists a unique solution to~\eqref{eq:problem_stat_2}, which we
  denote by $ \xnet_1^*=\ones_{n_1} \otimes x^*_1$ and $
  \xnet_2^*=\ones_{n_2} \otimes x^*_2 $.  By
  Proposition~\ref{prop:equiv-F}(i), there exists $
  (\xnet^*_1,\znet^*_1,\xnet^*_2,\znet^*_2) \in
  \equilibria{\directedNash}\cap W_{\znet^0_1,\znet^0_2} $.
  Consider the function $ \map{V}{(\real^{d_1n_1})^2\times
    (\real^{d_2n_2})^2}{\real_{\geq0}} $,
  \begin{align*}
    \Lyap(\xnet_1&,\znet_1,\xnet_2,\znet_2)
    \\
    &=\frac{1}{2}(\xnet_1-\xnet^*_1)^T(\xnet_1-\xnet^*_1)+
    \frac{1}{2}(\xnet_2-\xnet^*_2)^T(\xnet_2-\xnet^*_2)
    \\
    & \quad +
    \frac{1}{2}(\ynet_{(\xnet_1,\znet_1)}-\ynet_{(\xnet^*_1,\znet^*_1)})^T
    (\ynet_{(\xnet_1,\znet_1)}-\ynet_{(\xnet^*_1,\znet^*_1)}),
    \\
    & \quad +
    \frac{1}{2}(\ynet_{(\xnet_2,\znet_2)}-\ynet_{(\xnet^*_2,\znet^*_2)})^T
    (\ynet_{(\xnet_2,\znet_2)}-\ynet_{(\xnet^*_2,\znet^*_2)}),
  \end{align*}
  where $ \ynet_{(\xnet_\ell,\znet_\ell)}=\beta\xnet_\ell+\znet_\ell
  $, $ \ell \in \{1,2\} $, and $ \beta \in \realpositive $ satisfies $
  \beta^2-\alpha \beta+2=0 $.  This function is quadratic, hence
  smooth. Next, we consider $\xi =
  \Lie_{\directedNash}V(\xnet_1,\znet_1,\xnet_2,\znet_2)$ given by
  \begin{align*}
    \xi & =(-\alpha \Blap_1\xnet_1-\Blap_1 \znet_1 +\gradient
    \ULift(\xnet_1,\xnet_2),
    \Blap_1\xnet_1,-\alpha \Blap_2\xnet_2\\
    & \quad -\Blap_2 \znet_2 -\gradient \ULift(\xnet_1,\xnet_2),
    \Blap_2\xnet_2) \cdot \nabla V(\xnet_1,\znet_1,\xnet_2,\znet_2).
  \end{align*}
  After some manipulation, one can show that
  \begin{align*}
    \xi &= \sum_{\ell=1}^2\frac{1}{2} (\xnet_\ell-\xnet^*_\ell,
    \ynet_{(\xnet_\ell,\znet_\ell)}-\ynet_{(\xnet^*_\ell,\znet^*_\ell)})^T
    A_\ell ( \xnet_l, \ynet_{(\xnet_\ell,\znet_\ell)} )
    \\
    &+\sum_{\ell=1}^2\frac{1}{2}( \xnet_\ell^T,
    \ynet_{(\xnet_\ell,\znet_\ell)}^T ) A_\ell^T (
    \xnet_\ell-\xnet^*_\ell,
    \ynet_{(\xnet_\ell,\znet_\ell)}-\ynet_{(\xnet^*_\ell,\znet^*_\ell)})\\
    &+\sum_{\ell=1}^2(-1)^{j-1}(\xnet_\ell-\xnet^*_\ell)^T\gradient_{\xnet_\ell}
    \ULift(\xnet_1,\xnet_2)\\
    &+\sum_{\ell=1}^2(-1)^{j-1}\beta(\ynet_{(\xnet_\ell,\znet_\ell)} -
    \ynet_{(\xnet^*_\ell,\znet^*_\ell)})^T\gradient_{\xnet_\ell}
    \ULift(\xnet_1,\xnet_2),
  \end{align*}
  where $ A_\ell $, $ \ell \in \{1,2\} $, is
  \[
  A_\ell= 
  \begin{pmatrix}
    -(\alpha-\beta) \Blap_\ell & - \Blap_\ell\\
    (-\beta(\alpha-\beta)+1)\Blap_\ell & -\beta \Blap_\ell
  \end{pmatrix}.
  \]
  This equation can be written as
  \begin{align*}
    \xi &=\sum_{\ell=1}^2\frac{1}{2} ( \xnet_\ell-\xnet^*_\ell,
    \ynet_{(\xnet_\ell,\znet_\ell)}-\ynet_{(\xnet^*_\ell,\znet^*_\ell)}
    )^T
    \\
    &\qquad \qquad Q_\ell ( \xnet_\ell-\xnet^*_\ell,
    \ynet_{(\xnet_\ell,\znet_\ell)}-\ynet_{(\xnet^*_\ell,\znet^*_\ell)}
    )
    \\
    & + \sum_{\ell=1}^2 (\xnet_\ell-\xnet^*_\ell,
    \ynet_{(\xnet_\ell,\znet_\ell)}-\ynet_{(\xnet^*_\ell,\znet^*_\ell)}
    )^T A_\ell ( \xnet^*_\ell, \ynet_{(\xnet^*_\ell,\znet^*_\ell)} )
    \\
    &+\sum_{\ell=1}^2(-1)^{j-1}(\xnet_\ell-\xnet^*_\ell)^T\gradient_{\xnet_\ell}
    \ULift(\xnet_1,\xnet_2)
    \\
    &+\sum_{\ell=1}^2(-1)^{j-1}
    \beta(\ynet_{(\xnet_\ell,\znet_\ell)}-\ynet_{(\xnet^*_\ell,\znet^*_\ell)})^T
    \gradient_{\xnet_\ell} \ULift(\xnet_1,\xnet_2),
  \end{align*}
  where $ Q_\ell $, $ \ell \in \{1,2\} $, is given by
  \begin{equation}\label{eq:Q_with_y}
   Q_\ell=(\Blap_\ell+\Blap_\ell^T)\otimes \begin{pmatrix} 
     -(\frac{\beta^2+2}{\beta}-\beta) & -1\\
     -1 & -\beta
   \end{pmatrix}.
 \end{equation}   
   Note that, we have
   \begin{align*}
     A_1( \xnet^*_1, \ynet_{(\xnet^*_1,\znet^*_1)} )&=-(
     \Blap_1\ynet_{(\xnet^*_1,\znet^*_1)},
     \beta\Blap_1\ynet_{(\xnet^*_1,\znet^*_1)} )\\
     &= -(\gradient_{\xnet_1} \ULift(\xnet^*_1,\xnet_2),
     \beta\gradient_{\xnet_1} \ULift(\xnet^*_1,\xnet_2)),\\
     A_2( \xnet^*_2, \ynet_{(\xnet^*_2,\znet^*_2)} )&=-(
     \Blap_2\ynet_{(\xnet^*_2,\znet^*_2)},
     \beta\Blap_2\ynet_{(\xnet^*_2,\znet^*_2)} )\\
     &=(\gradient_{\xnet_2} \ULift(\xnet_1,\xnet^*_2),
     \beta\gradient_{\xnet_2} \ULift(\xnet_1,\xnet^*_2)).
  \end{align*}  
  Thus, after substituting for $ \ynet_{(\xnet_\ell,\znet_\ell)} $, we
  have
  \begin{align}\label{eq:z-issue}
    & \xi =\sum_{\ell=1}^2\frac{1}{2} (\xnet_\ell-\xnet^*_\ell,
    \znet_\ell-\znet^*_\ell )^T \tilde{Q}_\ell
    (\xnet_\ell-\xnet^*_\ell, \znet_\ell-\znet^*_\ell) \nonumber
    \\
    &\; +(1+\beta^2)(\xnet_1-\xnet^*_1)^T(\gradient_{\xnet_1}
    \ULift(\xnet_1, \xnet_2)-\gradient_{\xnet_1}
    \ULift(\xnet^*_1,\xnet_2)) \nonumber
    \\
    &\; -(1+\beta^2)(\xnet_2-\xnet^*_2)^T(\gradient_{\xnet_2}
    \ULift(\xnet_1, \xnet_2)-\gradient_{\xnet_2}
    \ULift(\xnet_1,\xnet^*_2)) \nonumber
    \\
    &\; +\beta(\znet_1-\znet^*_1)^T(\gradient_{\xnet_1}
    \ULift(\xnet_1,\xnet_2) - \gradient_{\xnet_1}
    \ULift(\xnet^*_1,\xnet_2)) \nonumber
    \\
    &\; -\beta(\znet_2-\znet^*_2)^T(\gradient_{\xnet_2}
    \ULift(\xnet_1,\xnet_2) - \gradient_{\xnet_2}
    \ULift(\xnet_1,\xnet^*_2)),
  \end{align}
  where 
  \begin{align*}
    \tilde{Q}_\ell=\begin{pmatrix}
      -\beta^3-(\frac{\beta^2+2}{\beta})-\beta & -(1+\beta^2)
      \\
      -(1+\beta^2) & -\beta
      \\
    \end{pmatrix}\otimes (\Blap_\ell+\Blap_\ell^T),
  \end{align*}
  for $ \ell \in \{1,2\} $. Each eigenvalue of $ \tilde{Q}_\ell $ is
  of the form
  \begin{align}\label{eq:eigenvalues-tildeQ}
    \tilde{\eta}_\ell = \lambda_\ell \frac{-(\beta^4+3\beta^2+2) \pm
      \sqrt{(\beta^4+3\beta^2+2)^2-4\beta^2}}{2\beta},
  \end{align}
  where $ \lambda_\ell $ is an eigenvalue of $ \Lap_\ell+\Lap_\ell^T
  $, $ \ell \in \{1,2\} $.  Using now Theorem~\ref{theorem:appendix}
  twice, one for $ (\xnet_1,\xnet_2) $, $ (\xnet_1^*,\xnet_2)$, and
  another one for $ (\xnet_1,\xnet_2) $, $ (\xnet_1,\xnet_2^*) $, we
  have
  \begin{align*}
    &(\xnet_1-\xnet^*_1)^T(\gradient_{\xnet_1} \ULift(\xnet_1,
    \xnet_2)-\gradient_{\xnet_1} \ULift(\xnet^*_1,\xnet_2))\leq\\
    &\quad -\frac{1}{K} \left(||\gradient_{\xnet_1} \ULift(\xnet_1,
      \xnet_2)-\gradient_{\xnet_1}
      \ULift(\xnet^*_1,\xnet_2)||^2\right),\\
    -&(\xnet_2-\xnet^*_2)^T(\gradient_{\xnet_2} \ULift(\xnet_1,
    \xnet_2)-\gradient_{\xnet_2} \ULift(\xnet_1,\xnet_2^*))\leq\\
    &\quad -\frac{1}{K} \left(||\gradient_{\xnet_2} \ULift(\xnet_1,
      \xnet_2)-\gradient_{\xnet_2}
      \ULift(\xnet_1,\xnet^*_2)||^2\right),
  \end{align*}
  where $ K \in \realpositive $ is the Lipschitz constant of $\nabla
  \ULift $.  We thus conclude that
  \begin{align*}
    &\xi \leq\sum_{\ell=1}^2\frac{1}{2} ( \xnet_\ell-\xnet^*_\ell,
    \znet_\ell-\znet^*_\ell)^T \tilde{Q}_\ell
    (\xnet_\ell-\xnet^*_\ell, \znet_\ell-\znet^*_\ell )
    \\
    &-\frac{(1+\beta^2)}{K} \big( ||\gradient_{\xnet_1}
    \ULift(\xnet_1, \xnet_2)-\gradient_{\xnet_1}
    \ULift(\xnet^*_1,\xnet_2)||^2 \\
    & +||\gradient_{\xnet_2} \ULift(\xnet_1,
    \xnet_2)-\gradient_{\xnet_2} \ULift(\xnet_1,\xnet^*_2)||^2 \big)
    \nonumber\\
    &+\beta(\znet_1-\znet^*_1)^T(\gradient_{\xnet_1}
    \ULift(\xnet_1,\xnet_2) - \gradient_{\xnet_1}
    \ULift(\xnet^*_1,\xnet_2))\\
    &-\beta(\znet_2-\znet^*_2)^T(\gradient_{\xnet_2}
    \ULift(\xnet_1,\xnet_2) - \gradient_{\xnet_2}
    \ULift(\xnet_1,\xnet^*_2)).\nonumber
  \end{align*}
  \newcounter{mytempeqncnt} One can write this inequality as displayed
  in~\eqref{eq:aux-display},
  \begin{figure*}[!htb]
    \normalsize
    \setcounter{mytempeqncnt}{\value{equation}}
    \begin{align}\label{eq:aux-display}
      \xi \leq \frac{1}{2}X^T
      \underbrace{\begin{pmatrix} 
          Q_{1,11} & Q_{1,12} & 0 & 0 &0 & 0 \\
          Q_{1,21} &Q_{1,22} & 0 & 0 & \beta \identity{n_1d_1} & 0 \\
          0 & 0 &Q_{2,11} & Q_{2,12} &0 & 0\\
          0 & 0 &Q_{2,21} &Q_{2,22} &0  & -\beta \identity{n_2d_2}\\
          0 & \beta \identity{n_1d_1} & 0 & 0 &
          -\tfrac{(1+\beta^2)}{K} \identity{n_1d_1}& 0
          \\ 
          0 & 0 & 0 &  -\beta \identity{n_2d_2} & 0 &
          -\tfrac{(1+\beta^2)}{K}\identity{n_2d_2} 
        \end{pmatrix}}_{\mathbf{Q}}
      X,  
    \end{align}
    \setcounter{equation}{\value{mytempeqncnt}}
    \hrulefill
  \end{figure*}    
  where
  \begin{align*}
    X=(&\xnet_1-\xnet^*_1,\znet_1-\znet^*_1,
    \xnet_2-\xnet^*_2,\znet_2-\znet^*_2,
    \\
    &\gradient_{\xnet_1} \ULift(\xnet_1, \xnet_2) -
    \gradient_{\xnet_1}\ULift(\xnet^*_1,\xnet_2),
    \\
    &\gradient_{\xnet_2} \ULift(\xnet_1, \xnet_2) -
    \gradient_{\xnet_2} \ULift(\xnet_1,\xnet^*_2)) .
  \end{align*}
  Since $ (\xnet_1, \znet_1,\xnet_2,\znet_2) \in
  W_{\znet^0_1,\znet^0_2} $, we have $(\ones_{n_\ell}^T \otimes
  \identity{d_\ell}) (\znet_\ell - \znet_\ell^*) =\zeros_{d_\ell}$, $
  \ell \in \{1,2\} $, and hence it is enough to establish
  that~$\mathbf{Q}$ is negative semidefinite on the subspace $\WW =
  \setdef{(v_1,v_2,v_3,v_4,v_5,v_6) \in
    (\real^{n_1d_1})^2\times(\real^{n_2d_2})^2\times
    \real^{n_1d_1}\times \real^{n_2d_2}}{(\ones_{n_1}^T \otimes
    \identity{d_1}) v_2 = \zeros_{d_1}, (\ones_{n_2}^T \otimes
    \identity{d_2}) v_4 =\zeros_{n_2}}$.  Using the fact that
  $-\tfrac{1}{K} (1+\beta^2) \identity{n_\ell d_\ell} $ is invertible,
  for $ \ell \in \{1,2\} $, we can express $\mathbf{Q}$ as
    \begin{align*}
    &\mathbf{Q}=\\
    &\ N
    \underbrace{\begin{pmatrix}
        \bar{Q}_1 & 0 & 0 &0\\
        0 & \bar{Q}_2 & 0 & 0
        \\
        0 & 0 & -\tfrac{1}{K} (1+\beta^2) \identity{n_1d_1} & 0\\
        0 & 0 & 0 & -\tfrac{1}{K} (1+\beta^2) \identity{n_2d_2}\\
      \end{pmatrix}}_{\mathbf{D}} N^T ,
  \end{align*}
  where $ \bar{Q}_\ell=\tilde{Q}_\ell+\frac{K\beta^2}{(1+\beta^2)}
  \left(\begin{matrix}
      0 & 0\\
      0 & \identity{n_\ell d_\ell}
    \end{matrix}\right)$, $ \ell \in\{1,2\} $, and 
  \begin{align*}
    N =
    \begin{pmatrix}
      \identity{n_1d_1} & 0 & 0 &0 &0 &0
      \\
      0 & \identity{n_1d_1} & 0 & 0 &-\frac{\beta K}{1+\beta^2}
      \identity{n_1d_1} &0 \\
      0 & 0 & \identity{n_2d_2} & 0 & 0 &0  \\
      0 & 0 & 0 & \identity{n_2d_2} & 0 & \frac{\beta K}{1+\beta^2}
      \identity{n_2d_2}
      \\
      0 & 0 & 0 & 0 & \identity{n_1d_1} & 0 \\
      0 & 0 & 0 & 0 & 0 & \identity{n_2d_2}\\
    \end{pmatrix}.
  \end{align*}
  Noting that $\WW$ is invariant under $N^T$ (i.e., $N^T \WW = \WW$),
  all we need to check is that the matrix $ \mathbf{D} $ is negative
  semidefinite on $\WW$. Clearly,
  \begin{align*}
    \left(
      \begin{smallmatrix}
        -\tfrac{1}{K} (1+\beta^2) \identity{n_1d_1} & 0
        \\
        0 & -\tfrac{1}{K} (1+\beta^2) \identity{n_2d_2}
      \end{smallmatrix} \right)    
  \end{align*}
  is negative definite. On the other hand, for $ \ell \in \{1,2\} $, on
  $(\real^{n_\ell d_\ell})^2$, $0$ is an eigenvalue of $\tilde{Q}_\ell$ with
  multiplicity $2d_\ell$ and eigenspace generated by vectors of the form
  $(\ones_{n_\ell} \otimes a,0)$ and $(0,\ones_{n_\ell} \otimes b)$, with
  $a,b \in \real^{d_\ell}$. However, on $ \setdef{(v_1,v_2) \in
    (\real^{n_\ell d_\ell})^2}{(\ones_{n_\ell}^T \otimes \identity{d_\ell}) v_2 =
    \zeros_{d_\ell}}$, $0$ is an eigenvalue of $\tilde{Q}_\ell$ with
  multiplicity $d_\ell$ and eigenspace generated by vectors of the form
  $(\ones_{n_\ell} \otimes a,0)$. Moreover, on $ \setdef{(v_1,v_2) \in
    (\real^{n_\ell d_\ell})^2}{(\ones_{n_\ell}^T \otimes \identity{d_\ell}) v_2 =
    \zeros_{d_\ell}}$, the eigenvalues of $\frac{K\beta^2}{(1+\beta^2)}
  \left(\begin{smallmatrix}
      0 & 0\\
      0 & \identity{n_\ell d_\ell}
    \end{smallmatrix}\right)
  $ are $\frac{K\beta^2}{(1+\beta^2)}$ with multiplicity $n_\ell d_\ell-d_\ell$
  and $0$ with multiplicity $n_\ell d_\ell$. Therefore, using Weyl's
  theorem~\citep[Theorem 4.3.7]{RAH-CRJ:85}, we deduce that the nonzero
  eigenvalues of the sum $ \bar{Q}_\ell $ are upper bounded by $
  \smallnonzeroeig(\tilde{Q}_\ell)+\frac{K\beta^2}{(1+\beta^2)} $. Thus,
  the eigenvalues of $\bar{Q}=\left(
    \begin{smallmatrix}
      \bar{Q}_1 & 0
      \\
      0 & \bar{Q}_2
    \end{smallmatrix} \right) $ are upper bounded by $
  \min_{\ell=1,2}\{\smallnonzeroeig(\bar{Q}_\ell)\}+\frac{K\beta^2}{(1+\beta^2)}
  $.  From~\eqref{eq:eigenvalues-tildeQ} and the definition of $h$
  in~\eqref{eq:h}, we conclude that the nonzero eigenvalues of~$
  \bar{Q} $ are upper bounded by $ h(\beta) $.  It remains to show
  that there exists $ \beta^* \in \realpositive $ with $ h(\beta^*)=0
  $ such that for all $ 0 < \beta < \beta^* $ we have $ h(\beta) < 0
  $.  For $ r>0 $ small enough, $h(r)<0$, since $ h(r) =
  -\frac{1}{2}\smallnonzeroeig^{\min}r+O(r^2)$. Furthermore, $
  \lim_{r\rightarrow \infty} h(r) =K > 0 $. Hence, using the Mean
  Value Theorem, we deduce the existence of $ \beta^* $.  Therefore we
  conclude that $ \Lie_{\directedNash}V
  (\xnet_1,\znet_1,\xnet_2,\znet_2) \leq 0$.  As a by-product, the
  trajectories of~\eqref{eq:Nash_dynamics-undirected} are bounded.
  Consequently, all assumptions of the LaSalle Invariance Principle,
  cf. Theorem~\ref{th:laSalle}, are satisfied. This result then implies
  that any trajectory of~\eqref{eq:Nash_dynamics-undirected} starting
  from an initial condition $
  (\xnet_1^0,\znet_1^0,\xnet_2^0,\znet_2^0) $ converges to the largest
  invariant set $ M $ in $ S_{\directedNash,\Lyap} \cap
  \WW_{\znet^0_1, \znet^0_2}$.  Note that if
  $(\xnet_1,\znet_1,\xnet_2,\znet_2) \in S_{\directedNash,\Lyap} \cap
  W_{\znet^0_1,\znet^0_2} $, then $ N^T X \in \ker (\bar{Q}) \times
  \{0\} $.  From the discussion above, we know $\ker (\bar{Q})$ is
  generated by vectors of the form $(\ones_{n_1} \otimes a_1,0,0,0)$,
  $(0,0,\ones_{n_2}\otimes a_2,0)$, $ a_\ell \in \real^{d_\ell} $, $ j
  \in \{1,2\} $, and hence 
  $\xnet_\ell = \xnet^*_\ell+\ones_{n_\ell}\otimes a_\ell $,
  $\znet_\ell = \znet_\ell^*$. Using the strict concavity-convexity,
  this then implies that $\xnet_\ell = \xnet^*_\ell $.
  Finally, for $(\xnet^*_1,\znet_1,\xnet^*_2,\znet_2) \in M$, using
  the positive invariance of $M$, one deduces that $ (\xnet^*_1,
  \znet_1,\xnet^*_2,\znet_2) \in \equilibria{\directedNash} $. \qed
\end{pf*}

\begin{remark}\longthmtitle{Assumptions on payoff
    function}\label{re:payoff-function-properties} {\rm Two
    observations are in order regarding the assumptions in
    Theorem~\ref{theorem:Nash-directed} on the payoff function. First,
    the assumption that the payoff function has a globally Lipschitz
    gradient is not too restrictive given that, since the state spaces
    are compact, standard boundedness conditions on the gradient imply
    the globally Lipschitz condition. Second, we restrict our
    attention to differentiable payoff functions because locally
    Lipschitz functions with globally Lipschitz generalized gradients
    are in fact differentiable,
    see~\cite[Proposition~A.1]{BG-JC:12-tac}.}
    \oprocend
\end{remark}

\begin{remark}\longthmtitle{Comparison with best-response
    dynamics}\label{remark:comp_best-response}
  {\rm Using the gradient flow has the advantage of avoiding the
    cumbersome computation of the best-response map. This, however,
    does not come for free.  There are concave-convex functions for
    which the (distributed) gradient flow dynamics, unlike the
    best-response dynamics, fails to converge to the saddle point,
    see~\citep{DF-FP:10} for an example.} \oprocend
\end{remark}


We finish this section with an example. 

\begin{example}\longthmtitle{Distributed adversarial
    selection of signal and noise power
    via~\eqref{eq:Nash_dynamics-directed}}\label{ex:example} 
  {\rm Recall the communication scenario described in
    Remark~\ref{remark:application}. Consider $ 5 $ channels,
    $\{\texttt{ch}_1,\texttt{ch}_2,\texttt{ch}_3,\texttt{ch}_4,\texttt{ch}_5\} $, for which the network $ \network_1 $
    has decided that $\{\texttt{ch}_1, \texttt{ch}_3\} $ have signal power $ x_1 $ and $
    \{\texttt{ch}_2,\texttt{ch}_4\} $ have signal power $ x_2 $. Channel $ \texttt{ch}_5 $ has its
    signal power determined to satisfy the budget constraint $ P\in
    \realpositive $, i.e., $P-2 x_1 -2x_2$. Similarly, the network $
    \network_2 $ has decided that $ \texttt{ch}_1 $ has noise power $y_1$,
    $\{\texttt{ch}_2,\texttt{ch}_3,\texttt{ch}_4\}$ have noise power $y_2$, and $ \texttt{ch}_5 $ has noise
    power $C-y_1 -3 y_2$ to meet the budget constraint $ C\in
    \realpositive $.  We let $ \xnet=(x^1,x^2,x^3,x^4,x^5) $ and $
    \ynet=(y^1,y^2,y^3,y^4,y^5) $, where $ x^i=(x^i_1,x^i_2)\in
    [0,P]^2 $ and $ y^i = (y^i_1,y^i_2)\in [0,C]^2 $, for each $ i\in
    \{1,\ldots,5\} $.  

    The networks $ \network_1 $ and $ \network_2 $, which are
    weight-balanced and strongly connected, and the engagement
    topology $\networkeng$ are shown in Figure~\ref{fig:ex_network}.
    \begin{figure}[htb!]
      \centering
      \includegraphics[width=.625\linewidth]{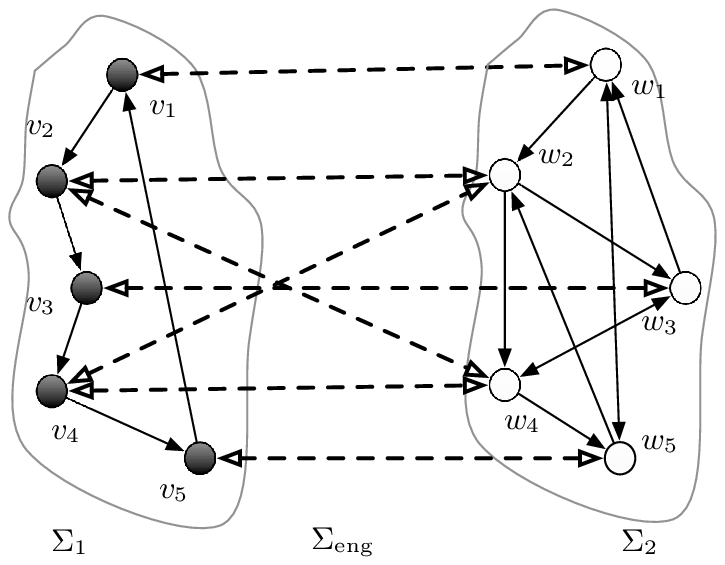}
      \caption{Networks $ \network_1 $, $ \network_2 $ and $
        \networkeng $ for the case study of Example~\ref{ex:example}.
        Edges which correspond to $ \networkeng $ are dashed. For 
        $ i\in \{1,\ldots, 5\}$, agents $v_i$ and $w_i$ are placed in
        channel $ \texttt{ch}_i $. }\label{fig:ex_network}
    \end{figure}
    Note that, according to this topology, each agent can observe the
    power employed by its adversary in its channel and, additionally,
    the agents in channel $ 2 $ can obtain information about the
    estimates of the opponent in channel $ 4 $ and vice versa.
    The payoff functions of the agents are given in
    Remark~\ref{remark:application}, where for simplicity we take $
    \sigma_i=\sigma_1 $, for $ i\in \{1,3, 5\} $, and $
    \sigma_i=\sigma_2 $, for $ i\in \{2,4\} $, with $
    \sigma_1,\sigma_2 \in \realpositive $.

    This example fits into the approach described 
    in Section~\ref{sec:reformulation} by considering the following 
    extended payoff functions:   
    \begin{align*}
      \fLift^1_1(x^1,\ynet) =&\log(1+\frac{\beta
        x^1_1}{\sigma_1+y^1_1}),
      \\
      \fLift^2_1(x^2,\ynet) =&\frac{1}{3}\log(1+\frac{\beta
        x^2_2}{\sigma_2+y^4_2}) +\frac{2}{3}\log(1+\frac{\beta
        x^2_2}{\sigma_2+y^2_2}),\\
      \fLift^3_1(x^3,\ynet) =&\log(1+\frac{\beta
        x^3_1}{\sigma_1+y^3_2}),\\
      \fLift^4_1(x^4,\ynet) =&\frac{1}{3}\log(1+\frac{\beta
        x^4_2}{\sigma_2+y^2_2})
      +\frac{2}{3}\log(1+\frac{\beta x^4_2}{\sigma_2+y^4_2}),\\
      \fLift^5_1(x^5,\ynet) =& \log \Big(1+ \frac{\beta (P-2x^5_1 -
        2x^5_2)}{\sigma_1+C-y^5_1 - 3 y^5_2} \Big),
      \\
      \fLift^1_2(\xnet,y^1)=&\fLift^1_1(x^1,\ynet), \quad
      \fLift^3_2(\xnet,y^3)=\fLift^3_1(x^3,\ynet),
      \\
      \fLift^2_2(\xnet,y^2) =&\frac{2}{3}\log(1+\frac{\beta
        x^2_2}{\sigma_2+y^2_2})
      +\frac{1}{3}\log(1+\frac{\beta x^4_2}{\sigma_2+y^2_2}),\\
      \fLift^4_2(\xnet,y^4) =&\frac{1}{3}\log(1+\frac{\beta
        x^2_2}{\sigma_2+y^4_2})
      +\frac{2}{3}\log(1+\frac{\beta x^4_2}{\sigma_2+y^4_2}),\\
      \fLift^5_2(\xnet,y^5)=&\fLift^5_1(x^5,\ynet).
    \end{align*}        
    Note that these functions are strictly concave and thus the
    zero-sum game defined has a unique saddle point on the set $
    [0,P]^2\times[0,C]^2 $.  These functions
    satisfy~\eqref{eq:properties}  and $ \ULift_1 = \ULift_2 $.
    %
    \begin{figure*}[htb!]
      \centering
      \subfigure[]{\includegraphics[width=.28\linewidth]{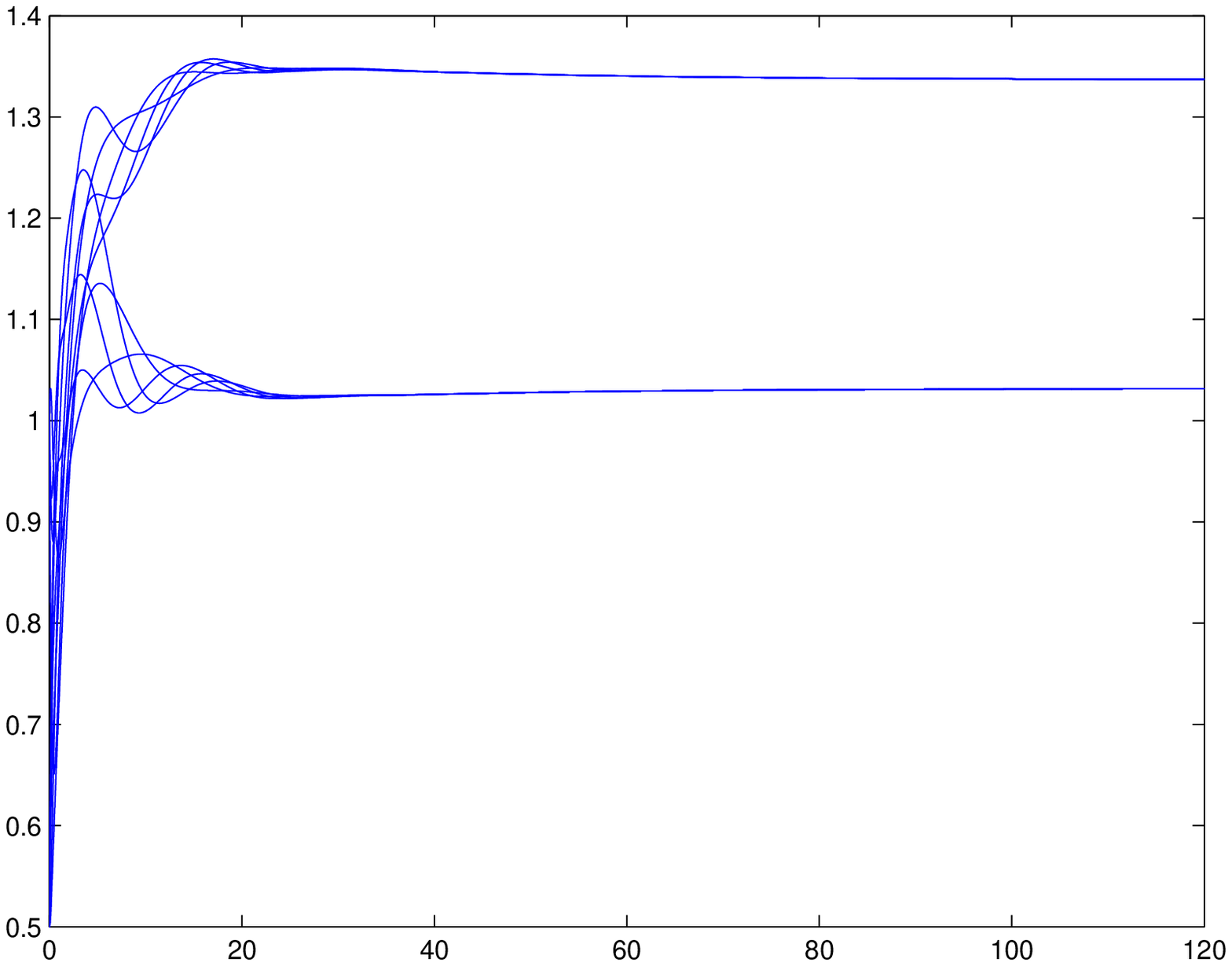}}
      \subfigure[]{\includegraphics[width=.28\linewidth]{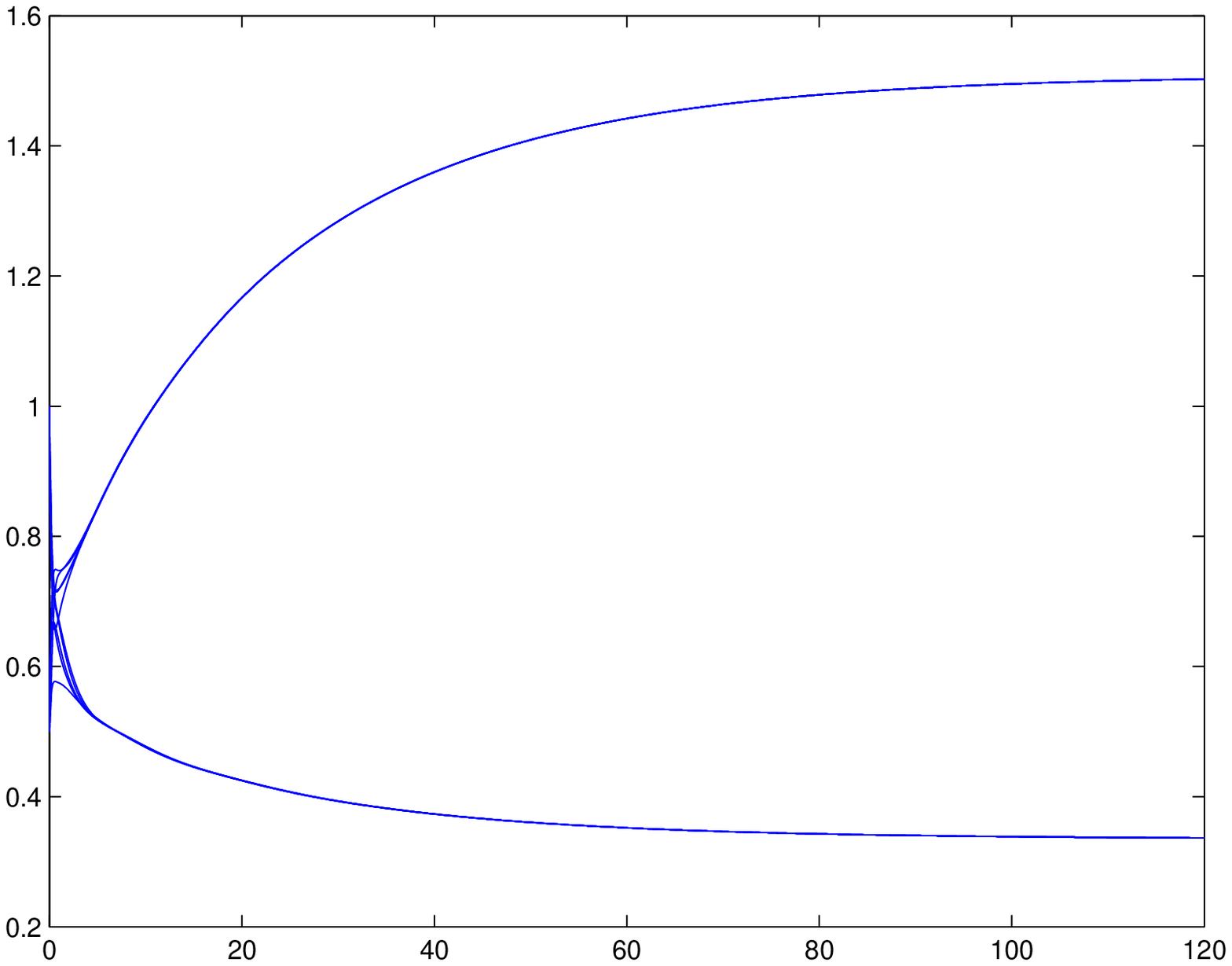}}
      \subfigure[]{\includegraphics[width=.28\linewidth]{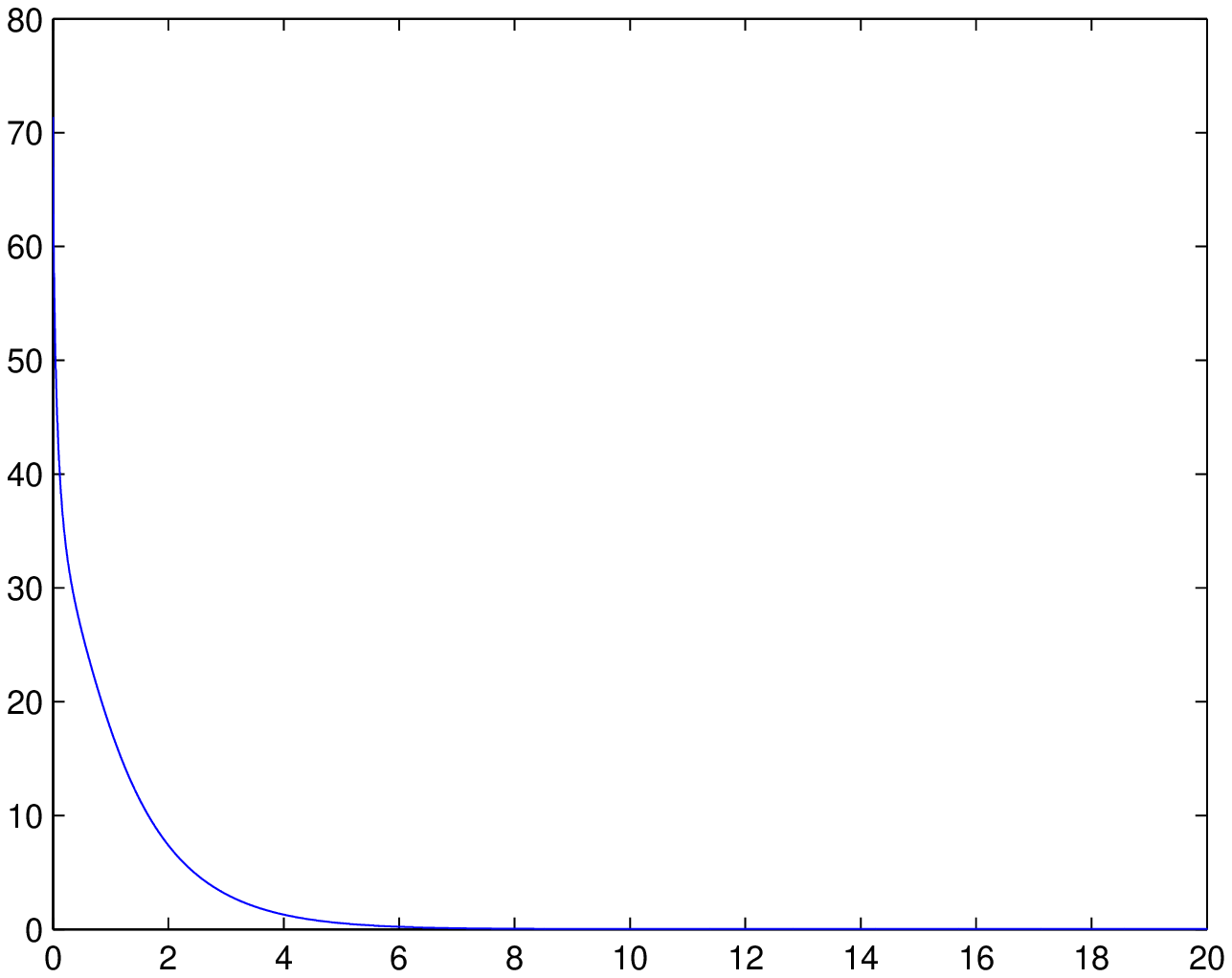}}
      \caption{Execution of~\eqref{eq:Nash_dynamics-directed} over the
        networked strategic scenario described in
        Example~\ref{ex:example}, with $ \beta=8 $, $ \sigma_1=1 $, $
        \sigma_2=4 $, $ P=6 $, and $ C=4 $.  (a) and (b) show the
        evolution of the agent's estimates of the state of networks $
        \network_1 $ and $ \network_2 $, respectively, and (c) shows
        the value of the Lyapunov function.  Here, $ \alpha =3 $
        in~\eqref{eq:Nash_dynamics-directed} and initially, $
        \xnet^0=((1,0.5),(0.5,1),(0.5,0.5),(0.5,1),(0.5,1))^T $, $
        \znet^0_1=\zeros_{10} $, $
        \ynet^0=((1,0.5),(0.5,1),(0.5,1),(0.5,0.5),(1,0.5))^T $ and $
        \znet^0_2=\zeros_{10} $. The equilibrium $
        (\xnet^*,\znet^*_1,\ynet^*,\znet^*_2) $ is $ \xnet^*=
        (1.3371,1.0315)^T\otimes\ones_5$, $ \ynet^* =
        (1.5027,0.3366)^T\otimes \ones_5 $, $ \znet^*_1 =
        (0.7508,0.5084,0.1447,0.5084,
        0.1447,-0.1271,-0.5201,-0.1271,-0.5201,-0.7626)^T$ and $
        \znet^*_2 =
        (0.1079,-0.0987,-0.0002,0.2237,0.0358,0.2875,-0.0360,
        0.0087,-0.1076,-0.4213)$.}\label{fig:sim}
    \end{figure*}
    Figure~\ref{fig:sim} shows the convergence of the
    dynamics~\eqref{eq:Nash_dynamics-directed} to the Nash equilibrium
    of the resulting $ 2 $-network zero-sum game.}  \oprocend
\end{example}

\myclearpage

\section{Conclusions and future work}\label{sec:conclusions}

We have considered a class of strategic scenarios in which two
networks of agents are involved in a zero-sum game. The networks aim
to either maximize or minimize a common objective function.
Individual agents collaborate with neighbors in their respective
network and have partial knowledge of the state of the agents in the
other one.  We have introduced two aggregate objective functions, one
per network, identified a family of points with a special saddle
property for this pair of functions, and established their
correspondence between the Nash equilibria of the overall game.  When
the individual networks are undirected, we have proposed a distributed
saddle-point dynamics that is implementable by each network via local
interactions. We have shown that, for a class of strictly
concave-convex and locally Lipschitz objective functions, the proposed
dynamics is guaranteed to converge to the Nash equilibrium.  We have
also shown that this saddle-point dynamics fails to converge for
directed networks, even when they are strongly connected and
weight-balanced.  Motivated by this fact, we have introduced a
generalization that incorporates a design parameter. We have shown
that this dynamics converges to the Nash equilibrium for strictly
concave-convex and differentiable objective functions with globally
Lipschitz gradients for appropriate parameter choices.  An interesting
venue of research is determining whether it is always possible to
choose the extensions of the individual payoff functions in such a way
that the lifted objective functions coincide.
Future work will also include relaxing the assumptions of strict
concavity-convexity and differentiability of the payoff functions, and
the globally Lipschitz condition on their gradients, extending our
results to dynamic interaction topologies and non-zero sum games, and
exploring the application to various areas, including collaborative
resource allocation in the presence of adversaries, strategic social
networks, collective bargaining, and collaborative pursuit-evasion.


\myclearpage

\myclearpage
\appendix

\section{Appendix}\label{sec:appendix}

The following result can be understood as a generalization of the
characterization of cocoercivity of concave
functions~\citep[Lemma~6.7]{EGG-NVY:96}.

\begin{theorem}\longthmtitle{Concave-convex differentiable functions
    with globally Lipschitz gradients}\label{theorem:appendix} Let $
  \map{f}{\real^{d_1}\times\real^{d_2}}{\real} $ be a concave-convex
  differentiable function with globally Lipschitz gradient (with
  Lipschitz constant $ K \in \realpositive $).  For $(x,y), (x',y')\in
  \real^{d_1}\times\real^{d_2} $,
  \begin{align*}
    (x& -x')^T (\gradient_xf(x,y)-\gradient_xf(x',y'))
    \\
    &\quad +(y-y')^T(\gradient_yf(x',y')-\gradient_yf(x,y))
    \\
    \leq & - \frac{1}{2K}
    \Big(\|\gradient_xf(x,y')-\gradient_xf(x',y')\|^2
    \\
    & +\|\gradient_yf(x',y)-\gradient_yf(x',y')\|^2
    \\
    &+\|\gradient_xf(x',y)-\gradient_xf(x,y)\|^2
    \\
    & + \| \gradient_yf(x,y')-\gradient_yf(x,y)\|^2 \Big) .
  \end{align*}
\end{theorem}

\begin{pf*}{Proof.}
  We start by noting that, for a concave function
  $\map{j}{\real^d}{\real}$ with globally Lipschitz gradient, the
  following inequality holds, see~\citep[Equation~6.64]{EGG-NVY:96},
  \begin{equation}\label{eq:G_property}
    j(x) \leq j^*-\frac{1}{2M}||\gradient j(x)||^2,
  \end{equation}
  where $ j^*=\sup_{x\in \real^{d}} j(x) $ and $M$ is the Lipschitz
  constant of $\nabla j$. Given $(x',y') \in \real^{d_1}\times
  \real^{d_2}$, define the map
  $\map{\tilde{f}}{\real^{d_1}\times\real^{d_2}}{\real} $ by
  \begin{align*}
    \tilde{f}(x,y)&=f(x,y)-f(x',y)-(x-x')^T\gradient_x f(x',y)
    \\
    & \quad +f(x,y)-f(x,y')+(y-y')^T\gradient_y f(x,y') .
  \end{align*}
  Since the gradient of $f$ is Lipschitz, the function $ \tilde{f} $
  is differentiable almost everywhere.  Thus, almost everywhere, we
  have
  \begin{align*}
    \gradient_x\tilde{f}(x,y) &= \gradient_xf(x,y)-\gradient_xf(x',y)
    + \gradient_xf(x,y)
    \\
    & \quad -
    \gradient_xf(x,y')-(y-y')^T\gradient_x\gradient_yf(x,y'),
    \\
    \gradient_y\tilde{f}(x,y) &= \gradient_yf(x,y) -
    \gradient_yf(x',y) + \gradient_yf(x,y)
    \\
    & \quad- \gradient_yf(x,y')-(x-x')^T\gradient_y\gradient_xf(x',y).
  \end{align*}
  In particular, note that $\gradient_x\tilde{f}(x',y') =
  \gradient_y\tilde{f}(x',y')=0$. Since $x \mapsto \tilde{f}(x,y')$
  and $y \mapsto \tilde{f}(x',y)$ are concave and convex functions,
  respectively, we can use~\eqref{eq:G_property} to deduce
  \begin{subequations}\label{eq:tilde_f_ineq}
    \begin{align}
      \tilde{f}(x,y') &\leq
      -\frac{1}{2K}||\gradient_xf(x,y')-\gradient_xf(x',y')||^2,
      \\
      -\tilde{f}(x',y) & \leq
      -\frac{1}{2K}||\gradient_yf(x',y)-\gradient_yf(x',y')||^2,
    \end{align}    
  \end{subequations}
  where we have used the fact that $\sup_{x\in \real^{d_1}}
  \tilde{f}(x,y') = \inf_{y \in \real^{d_2}}
  \tilde{f}(x',y)=\tilde{f}(x',y')=0 $.  Next, by definition of~$
  \tilde{f} $,
  \begin{align*}
    \tilde{f}(x,y')&=f(x,y')-f(x',y')-(x-x')^T\gradient_xf(x',y'),
    \\
    \tilde{f}(x',y)&=f(x',y)-f(x',y')-(y-y')^T\gradient_yf(x',y').
  \end{align*}
  Using~\eqref{eq:tilde_f_ineq}, we deduce that
  \begin{align}\label{eq:one_before}
    & f(x,y') -f(x',y) \nonumber
    \\
    & -(x-x')^T\gradient_xf(x',y') + (y-y')^T\gradient_yf(x',y')
    \nonumber
    \\
    & \leq -\frac{1}{2K} \big(\|\gradient_xf(x,y') -
    \gradient_xf(x',y')\|^2 \nonumber
    \\
    & \hspace*{2cm} +\|\gradient_yf(x',y)-\gradient_yf(x',y')\|^2
    \big) .
  \end{align}
  The claim now follows by adding together~\eqref{eq:one_before} and
  the inequality that results by interchanging $ (x,y) $ and $
  (x',y')$ in~\eqref{eq:one_before}. \qed
\end{pf*}


\begin{thebibliography}{31}
\providecommand{\natexlab}[1]{#1}
\providecommand{\url}[1]{\texttt{#1}}
\expandafter\ifx\csname urlstyle\endcsname\relax
  \providecommand{\doi}[1]{doi: #1}\else
  \providecommand{\doi}{doi: \begingroup \urlstyle{rm}\Url}\fi

\bibitem[Arrow et~al.(1951)Arrow, Hurwitz, and Uzawa]{KA-LH:51}
K.~Arrow, L~Hurwitz, and H.~Uzawa.
\newblock \emph{A Gradient Method for Approximating Saddle Points and
  Constrained Maxima}.
\newblock Rand Corporation, United States Army Air Forces, 1951.

\bibitem[Arrow et~al.(1958)Arrow, Hurwitz, and Uzawa]{KA-LH-HU:58}
K.~Arrow, L~Hurwitz, and H.~Uzawa.
\newblock \emph{Studies in Linear and Non-Linear Programming}.
\newblock Stanford University Press, Stanford, California, 1958.

\bibitem[Barron et~al.(2010)Barron, Goebel, and Jensen]{EN-RG-RRJ:10}
E.~N. Barron, R.~Goebel, and R.~R. Jensen.
\newblock Best response dynamics for continuous games.
\newblock \emph{Proceeding of the American Mathematical Society}, 138\penalty0
  (3):\penalty0 1069--1083, 2010.

\bibitem[Ba{\c s}ar and Olsder(1999)]{TB-GJO:99}
T.~Ba{\c s}ar and G.~J. Olsder.
\newblock \emph{Dynamic Noncooperative Game Theory}.
\newblock SIAM, 2 edition, 1999.
\newblock ISBN missing.

\bibitem[Boyd and Vandenberghe(2004)]{SB-LV:04}
S.~Boyd and L.~Vandenberghe.
\newblock \emph{Convex Optimization}.
\newblock Cambridge University Press, 2004.
\newblock ISBN 0521833787.

\bibitem[Bullo et~al.(2009)Bullo, Cort\'es, and Mart{\'\i}nez]{FB-JC-SM:08cor}
F.~Bullo, J.~Cort\'es, and S.~Mart{\'\i}nez.
\newblock \emph{Distributed Control of Robotic Networks}.
\newblock Applied Mathematics Series. Princeton University Press, 2009.
\newblock ISBN 978-0-691-14195-4.
\newblock Electronically available at http:/$\!$/coordinationbook.info.

\bibitem[Clarke(1983)]{FHC:83}
F.~H. Clarke.
\newblock \emph{Optimization and Nonsmooth Analysis}.
\newblock Canadian Mathematical Society Series of Monographs and Advanced
  Texts. Wiley, 1983.
\newblock ISBN 047187504X.

\bibitem[Cort{\'e}s(2008)]{JC:08-csm-yo}
J.~Cort{\'e}s.
\newblock Discontinuous dynamical systems - a tutorial on solutions, nonsmooth
  analysis, and stability.
\newblock \emph{{IEEE} Control Systems Magazine}, 28\penalty0 (3):\penalty0
  36--73, 2008.

\bibitem[Feijer and Paganini(2010)]{DF-FP:10}
D.~Feijer and F.~Paganini.
\newblock Stability of primal-dual gradient dynamics and applications to
  network optimization.
\newblock \emph{Automatica}, 46:\penalty0 1974--1981, 2010.

\bibitem[Frihauf et~al.(2012)Frihauf, Krstic, and Ba{\c s}ar]{PF-MK-TB:12}
P.~Frihauf, M.~Krstic, and T.~Ba{\c s}ar.
\newblock {N}ash equilibrium seeking in noncooperative games.
\newblock \emph{IEEE Transactions on Automatic Control}, 2012.
\newblock To appear.

\bibitem[Gharesifard and Cort\'es(2012{\natexlab{a}})]{BG-JC:11-acc}
B.~Gharesifard and J.~Cort\'es.
\newblock Distributed convergence to {N}ash equilibria by adversarial networks
  with undirected topologies.
\newblock In \emph{{A}merican {C}ontrol {C}onference}, pages 5881--5886,
  Montr\'eal, Canada, 2012{\natexlab{a}}.

\bibitem[Gharesifard and Cort\'es(2012{\natexlab{b}})]{BG-JC:12-cdc2}
B.~Gharesifard and J.~Cort\'es.
\newblock Distributed convergence to {N}ash equilibria by adversarial networks
  with directed topologies.
\newblock In \emph{{IEEE} Conf.\ on Decision and Control}, Maui, Hawaii,
  2012{\natexlab{b}}.
\newblock To appear.

\bibitem[Gharesifard and Cort\'es(2012{\natexlab{c}})]{BG-JC:12-tac}
B.~Gharesifard and J.~Cort\'es.
\newblock Continuous-time distributed convex optimization on directed graphs.
\newblock \emph{IEEE Transactions on Automatic Control}, 2012{\natexlab{c}}.
\newblock Conditionally accepted.

\bibitem[Golshtein and Tretyakov(1996)]{EGG-NVY:96}
E.~G. Golshtein and N.~V. Tretyakov.
\newblock \emph{Modified Lagrangians and Monotone Maps in Optimization}.
\newblock Wiley, New York, 1996.

\bibitem[Hofbauer and Sorin(2006)]{JH-SS:06}
J.~Hofbauer and S.~Sorin.
\newblock Best response dynamics for continuous zero-sum games.
\newblock \emph{Discrete and Continuous Dynamical Systems Ser. B}, 6\penalty0
  (1):\penalty0 215--224, 2006.

\bibitem[Horn and Johnson(1985)]{RAH-CRJ:85}
R.~A. Horn and C.~R. Johnson.
\newblock \emph{Matrix Analysis}.
\newblock Cambridge University Press, 1985.
\newblock ISBN 0521386322.

\bibitem[Johansson et~al.(2009)Johansson, Rabi, and Johansson]{BJ-MR-MJ:09}
B.~Johansson, M.~Rabi, and M.~Johansson.
\newblock A randomized incremental subgradient method for distributed
  optimization in networked systems.
\newblock \emph{SIAM Journal on Control and Optimization}, 20\penalty0
  (3):\penalty0 1157--1170, 2009.

\bibitem[Kim and Boyd(2008)]{JSK-SB:08}
J.~S. Kim and S.~Boyd.
\newblock A minimax theorem with applications to machine learning, signal
  processing, and finance.
\newblock \emph{SIAM Journal on Optimization}, 19\penalty0 (3):\penalty0
  1344--1367, 2008.

\bibitem[Li and Ba{\c s}ar(1987)]{SL-TB:87}
S.~Li and T.~Ba{\c s}ar.
\newblock Distributed algorithms for the computation of noncooperative
  equilibria.
\newblock \emph{Automatica}, 23\penalty0 (4):\penalty0 523--533, 1987.

\bibitem[Maistroskii(1977)]{DM:77}
D.~Maistroskii.
\newblock Gradient methods for finding saddle points.
\newblock \emph{Matekon}, 13:\penalty0 3--22, 1977.

\bibitem[Mesbahi and Egerstedt(2010)]{MM-ME:10}
M.~Mesbahi and M.~Egerstedt.
\newblock \emph{Graph Theoretic Methods in Multiagent Networks}.
\newblock Applied Mathematics Series. Princeton University Press, 2010.

\bibitem[Nedic and Ozdaglar(2009)]{AN-AO:09}
A.~Nedic and A.~Ozdaglar.
\newblock Distributed subgradient methods for multi-agent optimization.
\newblock \emph{IEEE Transactions on Automatic Control}, 54\penalty0
  (1):\penalty0 48--61, 2009.

\bibitem[Nedic and Ozdgalar(2009)]{AS-AO:09}
A.~Nedic and A.~Ozdgalar.
\newblock Subgradient methods for saddle-point problems.
\newblock \emph{Journal of Optimization Theory \& Applications}, 142\penalty0
  (1):\penalty0 205--228, 2009.

\bibitem[Olfati-Saber et~al.(2007)Olfati-Saber, Fax, and
  Murray]{ROS-JAF-RMM:07}
R.~Olfati-Saber, J.~A. Fax, and R.~M. Murray.
\newblock Consensus and cooperation in networked multi-agent systems.
\newblock \emph{Proceedings of the IEEE}, 95\penalty0 (1):\penalty0 215--233,
  2007.

\bibitem[Ren and Beard(2008)]{WR-RWB:08}
W.~Ren and R.~W. Beard.
\newblock \emph{Distributed Consensus in Multi-vehicle Cooperative Control}.
\newblock Communications and Control Engineering. Springer, 2008.
\newblock ISBN 978-1-84800-014-8.

\bibitem[Rockafellar(1997)]{RTR:97}
R.~T. Rockafellar.
\newblock \emph{Convex Analysis}.
\newblock Princeton Landmarks in Mathematics and Physics. Princeton University
  Press, Princeton, NJ, 1997.
\newblock ISBN 0-691-01586-4.
\newblock Reprint of 1970 edition.

\bibitem[Stankovic et~al.(2012)Stankovic, Johansson, and
  Stipanovic]{MSS-KHJ-DMS:12}
M.~S. Stankovic, K.~H. Johansson, and D.~M. Stipanovic.
\newblock Distributed seeking of {N}ash equilibria with applications to mobile
  sensor networks.
\newblock \emph{IEEE Transactions on Automatic Control}, 2012.
\newblock To appear.

\bibitem[Wan and Lemmon(2009)]{PW-MDL:09}
P.~Wan and M.~D. Lemmon.
\newblock Event-triggered distributed optimization in sensor networks.
\newblock In \emph{Symposium on Information Processing of Sensor Networks},
  pages 49--60, San Francisco, CA, 2009.

\bibitem[Wang and Elia(2010)]{JW-NE:10}
J.~Wang and N.~Elia.
\newblock Control approach to distributed optimization.
\newblock In \emph{Allerton Conf.\ on Communications, Control and Computing},
  pages 557--561, Monticello, IL, October 2010.

\bibitem[Wang and Elia(2011)]{JW-NE:11}
J.~Wang and N.~Elia.
\newblock A control perspective for centralized and distributed convex
  optimization.
\newblock In \emph{{IEEE} Conf.\ on Decision and Control}, pages 3800--3805,
  Orlando, Florida, 2011.

\bibitem[Zhu and Mart{\'\i}nez(2012)]{MZ-SM:12}
M.~Zhu and S.~Mart{\'\i}nez.
\newblock On distributed convex optimization under inequality and equality
  constraints.
\newblock \emph{IEEE Transactions on Automatic Control}, 57\penalty0
  (1):\penalty0 151--164, 2012.

\end{thebibliography}
\end{document}